\documentclass[english,12pt]{amsart}
\usepackage[utf8]{inputenc}
\usepackage{amssymb,stmaryrd}
\usepackage{amsfonts}
\usepackage{amscd}
\usepackage[all,cmtip]{xy}
\usepackage{color}
\usepackage{enumerate}
\usepackage{graphicx}
\usepackage{epstopdf} 
\usepackage[T1]{fontenc}
\RequirePackage[colorlinks=true, breaklinks=true, urlcolor= green, linkcolor= blue, citecolor= green, bookmarksopen=true,linktocpage=true,plainpages=false,pdfpagelabels]{hyperref}
 
\DeclareRobustCommand\longtwoheadrightarrow
    {\relbar\joinrel\twoheadrightarrow}

\newcommand{\CC}{\mathbb{C}}

\newcommand{\LL}{\mathbb{L}}
\newcommand{\NN}{\mathbb{N}}

\newcommand{\RR}{\mathbb{R}}

\newcommand{\ZZ}{\mathbb{Z}}

\newtheorem{thm}{Theorem}[section]
\newtheorem{cor}[thm]{Corollary}
\newtheorem{prop}[thm]{Proposition}
\newtheorem{lem}[thm]{Lemma}
\newtheorem{conj}[thm]{Conjecture}

\theoremstyle{definition}
\newtheorem{defn}[thm]{Definition}
\newtheorem{notn}[thm]{Notation}

\newtheorem{rem}[thm]{Remark}
\newtheorem{qu}[thm]{Question}
\newtheorem{ex}[thm]{Example}

\newcommand{\VF}{\operatorname{VF}}
\newcommand{\RF}{\operatorname{RF}}
\newcommand{\VG}{\operatorname{VG}}

\newcommand{\Th}{\mathrm{Th}}
\newcommand{\ac}{\mathrm{ac}}
\newcommand{\val}{\mathrm{val}}
\newcommand{\rtsp}{\mathrm{rtsp}}

\newcommand{\rad}{\operatorname{rad}}

\newcommand{\GL}{\operatorname{GL}}

\newcommand{\res}{\operatorname{res}}

\newcommand{\affgr}[1]{\bar{\operatorname{G}}_n^{#1}}
\newcommand{\vecgr}[1]{\operatorname{G}_n^{#1}}
\newcommand{\gr}[1]{\operatorname{G}_n^{#1}}
\newcommand{\OK}{\mathcal{O}_K}

\newcommand{\cC}{\mathcal{C}}
\newcommand{\Pt}{\operatorname{Pt}}

\newcommand{\cL}{\mathcal{L}}

\definecolor{better-green}{rgb}{0,.6,.1}
\definecolor{violet}{rgb}{.5,0,.5}

\long\def\changeR#1{{\color{red} #1}}
\long\def\changeG#1{{\color{better-green} #1}}

\newcommand{\private}[1]{}

\newcounter{dummy}

\begin{document}

\title[Motivic Vitushkin invariants]{Motivic Vitushkin invariants}
\author{Georges Comte}
\address{Univ.   Savoie Mont Blanc, CNRS, LAMA, 73000 Chamb\'ery, France}
\email{georges.comte@univ-smb.fr}
\urladdr{http://georgescomte.perso.math.cnrs.fr/}
\author{Immanuel Halupczok}
\address{Heinrich-Heine-Universit\"at D\"usseldorf,
Universit\"atsstr. 1, 40225 D\"usseldorf, Germany}
 \email{math@karimmi.de}
 \urladdr{https://immi.karimmi.de/}

\thanks{
\\
\\
The first author is supported by grant ANR DEFIGEO.
The 2nd author was partially supported by the research training group
\emph{GRK 2240: Algebro-Geometric Methods in Algebra, Arithmetic and Topology} and by the individual research grant No.~426488848, both funded by the Deutsche Forschungsgemeinschaft (DFG, German Research Foundation)}

\begin{abstract}
 Using motivic integration theory and the notion of riso-triviality, we introduce two new objects in the framework of definable nonarchimedean geometry: a convenient partial preorder $\preccurlyeq$ on the set of constructible motivic functions, extending the one considered in~\cite{CL.mot}, and an invariant $V_0$, nonarchimedean substitute of the number of connected components.
We then give several applications based on $\preccurlyeq$ and $V_0$: we obtain the existence of nonarchimedean substitutes of real measure geometric invariants $V_i$, called the Vitushkin variations, and we establish the nonarchi\-me\-dean counterpart of a real inequality involving $\preccurlyeq$, the metric entropy and our invariants $V_i$. We also prove the nonarchimedean Cauchy-Crofton formula for definable sets of dimension $d$, relating $V_0$ (and $V_d$) and the motivic measure in dimension $d$. 
\end{abstract}

\maketitle
\tableofcontents

\section*{Introduction}

In non-archimedean geometry, already 
one of the most basic geometric invariants from real geometry, namely the number of connected components of a set, becomes meaningless. The main aim of this article is to provide a substitute for this number: to each definable set $X$ in a (suitable) nonarchimedean context, we associate an invariant $V_0(X)$. The invariant $V_0(X)$ is not just a natural number, but contains some additional geometric information: It is a constant motivic function (which can be specialized to a natural number). A second aim is to
introduce a partial preorder $\preccurlyeq$ on the set of constructible motivic functions that, in the context of our motivic $V_0$, serves as a replacement for the usual order on the reals (the partial preorder from \cite{CL.mot} would not be suited for this purpose). 

To illustrate how $V_0$ and $\preccurlyeq$ allow new insights towards a consistent motivic integral geometry theory, we give several practical applications. Indeed, in the nonarchimedean context, from $V_0$ we define the counterpart of classical real integral geometry invariants, the Vitushkin variations $V_i$, and we prove the counterpart of celebrated formulas in real integral geometry, the Cauchy-Crofton formula, based on $V_0$, and the Ivanov formula for the metric entropy, based on the $V_i$'s and $\preccurlyeq$. 

In order to explain in detail why and how we develop this program, we begin in this introduction with a presentation of the real geometric context. 

In singularity theory, invariants defined through a deformation process and additive invariants (or maps)
are ubiquitous. Archetypical instances of such constructions are the Milnor fibre in singularity theory and the Lipschitz Killing curvatures in differential or integral geometry. 

The Milnor number of a 
polynomial function germ $f\colon (\CC^n,0)\to (\CC,0)$ having an isolated singularity at $0$ may be defined as (up to sign and shift by $1$) the Euler characteristic $\chi$ of the nearby fibre $\{f=\varepsilon\}\cap B(0,\eta)$, for $0<\varepsilon \ll \eta \ll 1$. Here, the nearby fibre $\{f=\varepsilon\}$
is viewed as a deformation of the singular fibre $\{f=0\}$. The Milnor number counts the cycles of real dimension $n-1$ in $\{f=\varepsilon\}\cap B(0,\eta)$ that vanish at the singular fibre ${f=0}$.

For $X$ a compact subset of $\RR^n$ definable in a given o-minimal structure expanding the real field, the Lipschitz-Killing curvatures of $X$
are defined through the deformation of $X$ by the family $(T_\varepsilon)_{\varepsilon >0}$ of its tubular neighbourhoods 
\[
T_\varepsilon(X) :=\{x\in \RR^n; 
\mathrm{dist}(x,X)\le \varepsilon\},
\]
and through the modified volume $\mathcal{V}_\varepsilon(X)$ of $T_\varepsilon(X)$, defined by 
\begin{equation}\label{eq.volume modifie}
\mathcal{V}_\varepsilon(X)=\int_{x\in T_\varepsilon(X)} \chi(B(x,\varepsilon)\cap 
X) \ \mathrm{d}x.
\end{equation}
Note that in case $X$ is smooth and compact, 
for $\varepsilon$  small enough, 
$\mathcal{V}_\varepsilon(X)=\mathrm{Vol}_n(T_\varepsilon(X))$, since 
$\chi(B(x,\varepsilon)\cap X)=1$ for any $\varepsilon>0$ small enough. 
Strikingly, it turns out that $\mathcal{V}_\varepsilon(X)$ is a polynomial in $\varepsilon$, i.e.,
\begin{equation}\label{eq.VolPol}
\mathcal{V}_\varepsilon(X)=\sum_{i=0}^n \alpha_i \Lambda_{n-i}(X) \varepsilon^i
\end{equation}
with universal coefficients $\alpha_i$ depending only on $i$, and coefficients $\Lambda_i(X)$ depending only on $X$. Those $\Lambda_i(X)$ are called the Lipschitz-Killing curvatures of $X$ (see for instance \cite{BerBro, Fu} for the definable case, and \cite{Wey} for the smooth case). 
On the other hand,  the $\Lambda_i(X)$'s are also defined by the following integral formulas 
\begin{equation}\label{eq.LipInt}
\forall i=0, \ldots, n, \ \ 
\Lambda_i(X)=\beta(i,n)\int_{\bar P\in  \affgr{n-i}}   \chi (X\cap \bar P )  \ \mathrm{d}\bar P,
\end{equation}
where $ \affgr{n-i}$ is the Grassmannian of $(n-i)$-dimensional
affine subspaces of $\RR^n$ and $\beta(i,n)$ is a universal constant depending only on $i$ and $n$ (see \cite{Com.singInvar} for a survey on this question).

A way to modify the Lipschitz-Killing curvature is to replace in
formula \eqref{eq.LipInt} the additive Euler characteristic $\chi$ by the sub-additive number of connected components $V_0$. In this way, one obtains a new sequence of invariants $V_0(X), \ldots, V_n(X)$ attached to $X$, called the Vitushin variations  (see for instance \cite{Iva, Vit55}, \cite[Chapter 3]{YoCo})
 \begin{equation}\label{eq.VarInt}
\forall i=0, \ldots, n, \ V_i(X)=\beta(i,n)\int_{\bar P\in \affgr{n-i} }   V_0 (X\cap \bar P )  \ \mathrm{d}\bar P.
\end{equation}
Note that directly from equality \eqref{eq.VarInt}, one obtains $V_n(X)=\mathrm{Vol}_n(X)$. It is a remarkable fact that when $\Lambda_i$ is replaced by $V_i$, (using $V_0$ instead of $\chi$ in \eqref{eq.LipInt}), 
the sequence  $V_0(X), \ldots, V_n(X)$ is still related to the volume of
$T_\varepsilon(X)$  by the notion of metric entropy $M(\varepsilon,X)$ of $X$, defined as the minimal number of balls of radius $\varepsilon>0$ needed to cover $X$. Indeed, a relation between the metric entropy and the variations of $X$, comparable to \eqref{eq.VolPol}
holds, but here, instead of having an equality, we only have a bound from above for $M(\varepsilon,X)$, due to the fact that $\chi$ has been replaced by the looser invariant~$V_0$:
\begin{equation}\label{eq.EntPol}
\forall \varepsilon>0, \  M(\varepsilon, X)\le C(n)\sum_{i=0}^n \frac{1}{\varepsilon^i}V_i(X),
\end{equation}
where $C(n)$ is a constant only depending on $n$ (see \cite{Iva}, \cite[Theorem 3.5]{YoCo}). 

To prove \eqref{eq.EntPol}, one introduces relative variations 
$V_i(X,B)$ of $X$, where $B$ is some ball. For this, one first defines 
$V_0(X,B)$ as the number of connected components of $X$ strictly lying in $B$ and then defines $V_1(X,B), \ldots, V_n(X,B)$ by 
\begin{equation}\label{eq.RelVar}
\forall i=0, \ldots, n, \ V_i(X,B)=\beta(i,n)\int_{\bar P\in \affgr{n-i} } V_0(X\cap \bar P, B)  \ \mathrm{d}\bar P.
\end{equation}
With this notation, inequality \eqref{eq.EntPol} is a consequence of the following one: 
\begin{equation}\label{eq.EntRelPol}
\sum_{i=0}^n V_i(X,B_1)\ge c(n),
\end{equation}
for a constant $c(n)$ depending only on $n$ and for any unit ball
$B_1$ centred at a point of $X$ 
 (see \cite[Theorem II.5.1]{Iva}, \cite[Theorem 3.4]{YoCo}, \cite[Theorem 2]{Zer} for a proof of \eqref{eq.EntRelPol}).
 
Inequality \eqref{eq.EntRelPol} says that the relative Vitushkin variations of $X$ in the unit ball $B_1$ cannot be too small all simultaneously. 
If $X$ has a connected component strictly included in $B_1$, then indeed, the sum of the $V_i(X,B_1)$'s is at least $1$ since $V_0(X,B_1)\ge 1$. Otherwise, the connected component of $X$ containing the centre of $B_1$ reaches the boundary of $B_1$, and this then yields that at least one relative variation $V_i(X,B_1)$, $i\ge 1$, is big. For instance, in the situation $\dim(X)=1$, a curve connecting the centre of $B_1$ and its boundary certainly lies in $X$, and $V_1(X,B_1)\ge V_1(R,B_1)=1$, where $R$ is a radius of $B_1$.

\medskip

In the definable nonarchimedean setting, the notion of deformation is also central. For instance, the notion of motivic Milnor fibre is defined using the rationality of a formal Poincaré series whose coefficients are based on a deformation of the singular fibre. The fact that such a Poincaré series is a rational function
shows that in this deformation process, some periodicity phenomenon and therefore some regularity is involved. Note the analogy to the polynomial formulas \eqref{eq.VolPol} and 
\eqref{eq.EntPol}, which also attest a moderation of the deformation process, viewed from the invariants $\Lambda_i(X)$ and $V_i(X)$. Moreover, in the rational form of the Poincaré series series, one can read some geometric data of the germ from which we start, namely (at least) its Milnor number.  

This general principle may be illustrated in the specific framework of this article, where the question of the existence of invariants related to the family of tubular neighbourhoods is addressed. 
For this, let us fix a complete equi-characteristic $0$ valued field $K$ with value group $ \ZZ$, and let us consider the definable family of tubular neighbourhoods 
 $ (T_r(X))_{r \in \NN}$ 
 of a definable subset $X$ of  the unit ball of $K^n$, namely 
 $$
 T_r(X):=\{y\in K^n ; \exists x\in X, \val(y-x)\ge r)\}.
 $$ 
By Theorem 14.4.1 of \cite{CL.mot}, the formal Poincaré series of the motivic volumes of this family is rational (see Section \ref{s.Motivic integrals} for notation):
\[
\sum_{r\in \NN} \mu_n(T_r(X))T^r \in \mathcal{C}(\{pt\})  \llbracket T \rrbracket_M, 
\]
where  $\mathcal{C}(\{pt\})  \llbracket T \rrbracket_M$ is the algebra of polynomials in $T$ with coefficients in $\mathcal Q(\{pt\})\otimes_{\ZZ[\LL-1]} A$
(see Notation \ref{not:motivic integration}) localized by the multiplicative set $M$ generated by polynomials 
$1-\LL^\alpha T^\beta$, $\alpha\in \ZZ, \beta\in \NN\setminus\{0\}$. 
In particular, considering a generator
$\dfrac{P_k(T)}{(1-\LL^{\alpha_k} T^{\beta_k})^k}$ of $\mathcal{C}_{\{pt\}}  \llbracket T \rrbracket_M$, its Taylor formula shows that it contributes to the series only for coefficients of order $r=p+m\beta_k$, $m\in \NN$, where $p$ runs over the degrees of the monomials appearing in $P_k$. Moreover, this contribution has the form $\LL^{m\alpha_k} Q_k(m) $ with $Q_k$ a polynomial with coefficients in $\mathcal Q({\{pt\}})\otimes_{\ZZ[\LL-1]} A$. 
In conclusion, there exists a finite number of 
functions of $r$,  
$R_1, \ldots, R_\ell$
such that   

\begin{equation}\label{eq.rationality}
\mu_n(T_r(X))\in \{ R_1(r), \ldots, R_\ell(r)
\},
\end{equation}
each $R_i$ being a linear combination
with coefficients in $\mathcal Q({\{pt\}})\otimes_{\ZZ[\LL-1]} A$ of elements of the form 
$$
C_k(\frac{r-p_{k,j}}{\beta_k})\LL^{\alpha_k\frac{r-p_{k,j}}{\beta_k}},
$$ 
where 
$C_k$ are polynomials, $p_{k,j}$ range in a finite subset of $\NN$, $\alpha_k\in \ZZ$, $ \beta_k\in \NN$, with the condition that $r=p_{k,j}\mod{\beta_k}$.

The functions $R_1, \ldots, R_\ell$ encode the arithmetico-logical complexity of the set $X$, but it does not seem that one can easily extract from them some geometrical meaning, in contrast to the case of the Poincaré series defining the motivic Milnor fibre, and in contrast to the smooth or  convex real case, where $\ell=1$ and $R_1$ is a polynomial function of the radius of the tubular neighbourhood, from which the Lipschitz Killing curvatures appear as coefficients (see \eqref{eq.VolPol}). In the tubular deformation family no motivic version of the $\Lambda_i$'s may clearly be identified, except maybe in the trivial case where $X$ is a smooth variety of dimension $d$: In this case only the $d$-volume of $X$ appears since one has $\mu_n(T_r(X))=\mu_d(X)\LL^{-r(n-d)}$.
Note finally that since in the Hrushovski--Kazhdan version of motivic integration \cite{HK.motInt} a notion of motivic Euler-characteristic 
appears, one may think that defining with the latter a motivic modified volume modelled on formula \eqref{eq.volume modifie} could make it possible to relate it to
a motivic version of the $\Lambda_i$'s, also defined in this context by \eqref{eq.LipInt}, for instance by a formula comparable to \eqref{eq.VolPol}.
But so far all attempts in this direction failed.

Nevertheless, in spite of this apparent lack of geometric information in the deformation of $X$ via the family of tubular neighbourhoods,   
we show in Section \ref{s.Motivic Vitushkin invariants} that one can transfer to the nonarchimedean setting
the notion of Vitushkin variations, and in Section \ref{s.Motivic Vitushkin formula} that these motivic Vitushkin variations are this time deeply related to the motivic metric entropy.  
Indeed, one of the main results of this article is Theorem \ref{t.ent}, the nonarchimedean counterpart of inequality
\eqref{eq.EntPol}. 

By definition, the motivic metric entropy $M(X,r)$ of $X$
is the motivic measure (instead of the number in the real situation) of balls one needs to cover $X$ (see Definition \ref{d.Motivic entropy}). But it turns out in the nonarchimedean setting that $M(X,r)$ is also $\mu(T_r(X))$, and thus
the family $(M(X,r))_{r\in \NN}$ is in general expressed in the form given by \eqref{eq.rationality}, from which there is {\sl a priori} no obvious way to identify some convenient geometrical invariants attached to~$X$. 

The principal obstruction we have to overcome when one aims at finding nonarchimedean substitutes of the variations $V_i(X)$, is finding a sensible notion of (number of) connected components. For this, we use the notion of riso-triviality introduced in \cite{iBW.canStrat} by Bradley-Williams and the second author (see Section \ref{s.riso-triviality}). 
Intuitively, $V_0$ will be the $0$-dimensional motivic measure of the $0$-dimensional stratum of a (certain kind of) minimal stratification of $X$.
While one can think of a t-stratification in the sense of \cite{i.whit}, those may be inherently non-canonical; to obtain a well-defined $V_0$, we use the more abstract approach from \cite{iBW.canStrat} instead.

To define a nonarchimedean version of the relative variation $V_0(X, B)$ we also need a replacement for the notion of connected component of $X$ lying strictly in $B$. It turns out that as a replacement for this, we can use the condition that $X$ can be stratified in such a way that $B$ is disjoint from the $0$-dimensional stratum.
For the accurate definition of the nonarchimedean version $V_0$, see Definitions \ref{d.V0} and \ref{d.V_0 relative}.

Another technical difficulty to establish the nonarchimedean version of inequality \eqref{eq.EntPol} or \eqref{eq.EntRelPol}, is the following. 
Since $V_0(X)$ is defined as the $0$-dimensional motivic measure of some definable set in the valued field,
we have to introduce in Section \ref{s.preorder} a notion of partial preorder allowing us to compare for instance a class in $\mathcal Q_+({\{pt\}})$, the Grothendieck semiring of definable sets in the residue field sort, and 
a constant function with value in $\ZZ$ (typically $[x^2-2=0]$ and the constant function of $\mathcal{C}_+(\{pt\})$ equal to $1$ (see Section \ref{s.Motivic integrals} for the notation).  
This crucial need explicitly appears in this simple form in the proof of Theorem \ref{t.bound on sum of variations}. 
The partial preorder introduced in 
\cite[12.2]{CL.mot}, defined via the set of nonnegative constructible motivic functions (see Remark \ref{r.positivity in P+}), is in general not adapted to this goal. 
Our more flexible notion of partial preorder on the set of integrable constructible motivic functions, denoted $\preccurlyeq$ and extending the one of \cite[12.2]{CL.mot} is introduced in Definition \ref{d.<} and Lemma \ref{l.<+}, using motivic integration and definable surjections. 

With the invariants $V_i$ defined through the nonarchimedean invariant $V_0$ as  in the real case (see Definition \ref{d.Vi}), as well as their relative version $V_i(\cdot, B)$ in a ball $B$ (see Definition \ref{d.V_0 relative}), and with this partial preorder relation $\preccurlyeq$, 
we prove in Section \ref{s.Motivic Vitushkin formula}
 the nonarchimedean version of inequality \eqref{eq.EntRelPol} (see Theorem \ref{t.bound on sum of variations}) and finally 
the nonarchimedean version of inequality \eqref{eq.EntPol} (see Theorem \ref{t.ent}) relating the motivic entropy defined in Definition \ref{d.Motivic entropy} and the motivic Vitushkin variations $V_i$.

It also turns out (see Theorem \ref{t.Crofton}), in case $X$ is a definable subset of $K^n$ of dimension $d$, that $V_d(X)$ is the motivic 
$d$-dimensional volume of $X$ (up to some universal constant depending only on $d$ and $n$). In other words, the motivic Cauchy-Crofton formula involving 
$V_0$ as the nonarchimedean substitute of the real measure counting holds. This application of the construction of $V_0$ shows, at least from the geometric motivic measure theory point of view, that the invariant $V_0$ is a good substitute of the real number of points. 
Such a formula in the nonarchimedean context has been obtained in the local case by A. Forey in \cite{For.LocalCauchy}, on the model of the real local Cauchy-Crofton formula previously obtained in \cite{Co1} for the density by the first author. In this paper, we show the classical (global) motivic version of this formula, the proof of which requires a different approach than in the local case.  

The nonarchimedean Cauchy-Crofton formula, as well as the nonarchimedean version of inequality \eqref{eq.EntPol}, may be viewed as a step towards motivic integral geometry formulas in the nonarchimedean context, based on the $0$-degree motivic invariant $V_0$, such as classical kinematic formulas for homogeneous spaces (note that in the framework of characteristic zero local fields integral geometry tools have been  recently developed using $p$-adic integration in \cite{BuKuLe, KuLe}). In this direction, our preorder has a role to play, since the absence of a suitable nonarchimedean substitute of the Euler characteristic might make it necessary to replace equalities in classical real formulas involving the Euler characteristic by inequalities involving $V_0$. As explained above, this principle is here applied to the lack of a geometric nonarchimedean version of formula  \eqref{eq.VolPol}, turned into the geo\-metric nonarchimedean version of inequality \eqref{eq.EntPol}.

\section{Language, Theory, Structure, definable sets}\label{s.Language}

\begin{defn}
In the entire paper, $k$ is a field of characteristic $0$ and $K := k(\!(t)\!)$ is the field of formal Laurent series.
We write $\OK = k[\![t]\!]$ for the valuation ring of $K$, 
$\res\colon \OK \to k$ for the residue map,
$\ac\colon K \to k$ for the (natural) angular component map and $\val\colon K \to \ZZ \cup \{\infty\}$ for the valuation map, 
$\mathcal{M}_K=\{ x\in K; \val(x)>0 \}$, the maximal ideal.
We extend the valuation to $K^n$, by setting
$\val((a_1, \dots, a_n)) := \min_i \val(a_i)$.
\end{defn}

\begin{notn}\label{n.norm}
We sometimes also use a more intuitive, multiplicative notation for the valuation: \[
K \to \LL^\ZZ \cup \{0\},
a \mapsto |a| := \LL^{-\val(a)},
\]
and the multiplicative value group $\LL^\ZZ$ is equipped with the order satisfying $\LL^r \le \LL^s$ if and only if $r \le s$. Note that this notation suggests a natural map from $\LL^\ZZ$ into the set of constructible motivic functions $\mathcal{C}(\{pt\})$ on the point (see Notation \ref{not:motivic integration}).
\end{notn}

\begin{notn}
For $x \in K^n$ and $\lambda \in \LL^\ZZ$, we write 
$$
B(x,\lambda) = \{y \in K^n; |y - x| \le \lambda\}
$$
for the closed ball of radius $\lambda$ around $x$. By a \emph{ball}, we mean such a closed valuative ball. We denote 
the radius $\lambda$ of $B(x,\lambda) $ by $\rad B(x,\lambda)$. 
\end{notn}

For the entire paper, we also fix a language $\cL$ and consider $K$ as an $\cL$-structure. One possibility is to take $\cL$ to be the language $\cL_{\mathrm{DP}}$ of Denef--Pas.

\begin{defn}
The language $\cL_{\mathrm{DP}}$ of Denef--Pas consists of one sort for $K$ with the ring language, one sort for $k$ with the ring language, one sort for $\ZZ$ (the value group) with the language $\{0,+, <\}$ of ordered abelian groups, the valuation map and the angular component map.
\end{defn}

However, all our results also hold for various other language $\cL$. Here are the assumptions we really need: 

\begin{defn}\label{d.cL}
Let $\cL$ be an expansion of the language $\cL_{\mathrm{DP}}$ on $K$ with the following properties
\begin{enumerate}
 \item The theory $\Th_{\cL}(K)$ is $1$-h-minimal in the sense of \cite{iCR.hmin}.
 \item The value group and the residue field are orthogonal (i.e. every definable set $X \subseteq k^m \times \ZZ^n$ is a finite union of sets of the form $Y \times Z$, for $Y \subseteq k^m$ and $Z \subseteq \ZZ^n)$.
 \item The induced structure on $\ZZ$ is the pure ordered abelian group structure.
\end{enumerate}
\end{defn}

\begin{rem}
All conditions of Definition~\ref{d.cL} are preserved under adding constants from $K$, $k$ or $\ZZ$, so anything which we prove about $\emptyset$-definable sets and/or maps also holds for $A$-definable sets/maps, when $A$ is a subset of $K\cup k \cup \ZZ$. We will use this implicitly various times.
\end{rem}

Instead of imposing $K = k(\!(t)\!)$, one could have allowed $K$ to be a Henselian valued field with residue field $k$ and value group $\ZZ$. The reasons not to do that are that (a) if $K$ is not spherically complete, the definition of $V_0$ becomes more technical and (b) the case $K = k(\!(t)\!)$ implies the more general case, as follows.

First note that any such more general $K$ can be embedded into its maximal immediate extension, which is isomorphic to $k(\!(t)\!)$, and by \cite{iBW.sph}, the $\cL$-structure on $K$ can be extended to $k(\!(t)\!)$ in such a way that $k(\!(t)\!)$ is an elementary extension of $K$. For that reason, assuming $K = k(\!(t)\!)$ is not a restriction, provided that the results we are proving are first order statements. This is indeed the case for our results, in the following sense.

Our definition of $V_0$ (Definition~\ref{d.V0}) uses the notion of riso-triviality (Definition~\ref{d.risotriv}) which, as formulated, is not a first order condition. Moreover, if $K$ is not spherically complete, defined like this, the notion of riso-triviality would lose many of its good properties, making it useless for our purposes  (see Remark~\ref{r.riso.surj}). We do not know a direct way to fix this, but there is a simple indirect fix. By Theorem~\ref{t.defble}, there exists a formula expressing riso-triviality in every spherically complete field. We simply take that formula as the definition of riso-triviality in non-spherically complete fields, thus making that notion first order by definition.

With this definition of $V_0$, all our results (in particular Theorem~\ref{t.Crofton} and Theorems~\ref{t.bound on sum of variations} and \ref{t.ent}) also hold
in any $K \prec k(\!(t)\!)$. Indeed, those are equalities and inequalities between constructible motivic functions. Any such (in)equality is witnessed by finitely many formulas (e.g.\ defining some bijections or ingredients of other constructible motivic functions); the witnesses that work in $k(\!(t)\!)$ then also work in $K$.
\begin{rem}\label{r.k.st.e}
\cite[Proposition 2.6.12]{iCR.hmin} states that under the assumption of $1$-h-minimality, the $\mathrm{RV}$-sort (which, in our setting, is just $k \times \ZZ$) is stably embedded in a strong sense. Combining this with orthogonality of the value group and the residue field, one deduces that also the residue field is stably embedded in that strong sense.
One way to formulate this is the following: If $X \subset K^n \times k^m \times \ZZ^r$
and $Z \subseteq X \times k^{m'}$ are $\emptyset$-definable, for some $n, m, r, m' \in \NN$, then there exists an $m'' \in \NN$, an $\emptyset$-definable map $\rho\colon X \to k^{m''}$ and a $\emptyset$-definable set $\hat Z \subseteq k^{m''} \times k^{m'}$ such that for every $x \in X$, the fibre $Z_x$ is equal to the fibre $\hat Z_{\rho(x)}$ (to get from the version of stable embeddedness of \cite[Proposition 2.6.12]{iCR.hmin} to the one given here, one uses a standard compactness argument: First, one applies stable embeddedness to each fibre $Z_x$ individually, and then compactness shows that the resulting $\rho(x)$ and $\hat Z_{\rho(x)}$ can be taken to be definable uniformly in $x$).
\end{rem}

\section{Riso-triviality}\label{s.riso-triviality}

In this section we present the notion of riso-triviality defined and studied in \cite{iBW.canStrat}, a notion that we use in the rest of the paper. A set $X\subseteq K^n$ is said to be riso-trivial on some ball $B \subseteq K^n$ if it is roughly translation invariant on $B$, in some directions specified by a sub-vector space $V \subseteq k^n$. Note that many definitions and statements in \cite{iBW.canStrat} speak about riso-triviality of maps $\chi$ defined on $K^n$. If one is interested in a set $X \subseteq K^n$, one should take $\chi$ to be the indicator function of $X$.

\begin{defn}\label{d.riso}
A map $\Phi\colon X \to Y$ between sets $X, Y \subseteq K^n$ is called a \emph{risometry} if for every $x, x' \in X$ with $x \ne x'$, we have
$$\val((\Phi(x')-\Phi(x))-(x'-x))>\val(x'-x).$$
\end{defn}

\begin{rem}\label{r.riso.surj}
Clearly, risometries are always injective. In our setting, a risometry from a ball $B \subseteq K^n$ to itself is moreover always surjective. Indeed, since $K = k(\!(t)\!)$ is spherically complete (meaning that any descending chains of balls has non-empty intersection), surjectivity can be deduced using a version of the Banach Fixed Point Theorem in such fields (see \cite[Lemma~3.3.4]{iBW.canStrat}).
That surjectivity is crucial for the existence of $\rtsp_B$ in Definition~\ref{d.risotriv} below.
\end{rem}

\begin{defn}[Riso-triviality on a ball]
\label{d.risotriv}
Suppose that $X \subseteq K^n$ is an arbitary set and that $B \subseteq K^n$ is a ball.
\begin{enumerate}
 \item 
Given a sub-vector space $V$ of $k^n$, we
say that $X$ is \emph{$V$-riso-trivial on $B$} if there exists a lift $\tilde V \subseteq K^n$ of $V$ (i.e. satisfying $\res(\tilde V) = V$) and a risometry $\Phi\colon B \to B$ such that $\Phi^{-1}(X \cap B)$ is $\tilde V$-translation invariant, i.e.
for every pair of points $x,y\in B$ satisfying $y-x\in \tilde{V}$, we have 
$$ x \in \Phi^{-1}(X\cap B)\Longleftrightarrow y\in \Phi^{-1}(X\cap B).$$ 
\begin{figure}
    \centering
  \hskip1cm  
  \includegraphics[scale=0.55]{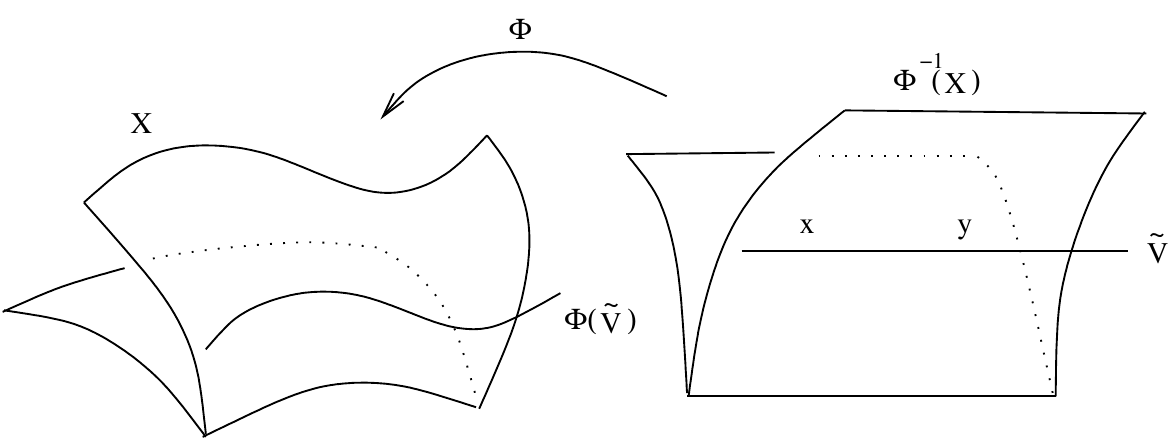}
  \vskip0cm
    \caption{$1$-riso-triviality on a ball}
    \label{f.risotriviality}
\end{figure}

\item
The \emph{riso-triviality space} of $X$ on $B$ is the maximal (with respect to inclusion) subspace $V \subseteq k^n$ such that $X$ is $V$-riso-trivial on $B$; we denote it by $\rtsp_B(X)$. Such a maximal subspace exists by \cite[Proposition~2.3.12]{iBW.canStrat}.
\item
We say that $X$ is \emph{$d$-riso-trivial} on $B$ if $\dim \rtsp_B(X) \ge d$, and we say that it is not riso-trivial if $X$ is not even $1$-riso-trivial. If it is clear which set $X$ we are talking about, we may also say that $B$ is $d$-riso-trivial, meaning that $X$ is $d$-riso-trivial on $B$.
\end{enumerate}

\end{defn}

Note that even though the above notions will be applied to definable sets $X$, we do not require the risometry $\Phi$ to be definable.

We will apply the above notions also for $X$ and $B$ living in an affine subspace $\bar P \subseteq K^n$, by identifying $\bar P$ with $K^{\dim \bar P}$.

Various good properties of these notions are proved in \cite{iBW.canStrat}. We now recall a few of them.

\begin{rem}
One easily sees that if a set $X \subseteq K^n$ is $V$-riso-trivial on a ball $B \subseteq K^n$, for some $V \subseteq k^n$, then $X$ is also $V'$-riso-trivial on any sub-ball $B' \subseteq B$ for any sub-vector space $V' \subseteq V$.
\end{rem}

\begin{rem}
In Definition~\ref{d.risotriv} (1), the choice of the lift $\tilde V$ does not matter, i.e. if there exists a lift and a risometry $\Phi$ witnessing $V$-riso-triviality, then for any other lift $\tilde V'$ of $V$, there exists a risometry $\Phi'$ such that $(\tilde V', \Phi')$ also witness $V$-riso-triviality.
\end{rem}

\begin{prop}[see {\cite[Lemma~3.3.14]{iBW.canStrat}}]\label{p.fib}
Let $X\subseteq K^n$ be a set, let $B\subseteq K^n$ be a ball, and set $\ell := n - \dim \rtsp_B(X)$.
Fix $\bar P$ in the Grassmannian manifold $\affgr{\ell}$ of $\ell$-dimensional affine subspaces of 
$K^n$, such that $\bar P \cap B \ne \emptyset$
and such that, if we denote by $P \subseteq K^n$ the direction of $\bar P$,
$\res(P) \cap \rtsp_B(X) = \{0\}$. Then $X \cap \bar P$ is not riso-trivial on $B \cap \bar P$.
\end{prop}

The idea of the proof of this proposition is that
if $X$ would be $V$-riso-trivial on $B \cap \bar P$, this would also imply $V$-riso-triviality on $B \cap \bar P'$ for all $\bar P'$ parallel to $\bar P$ and we could deduce $(V\oplus  \rtsp_B(X))$-riso-triviality of $X$ on $B$.

One of the central results of \cite{iBW.canStrat} is that $\rtsp_B(X)$ depends definably on $X$ and $B$. More precisely, we have the following statement. 

\begin{thm}[{\cite[Corollary~3.1.3]{iBW.canStrat}}]\label{t.defble}
Suppose we are given a $\emptyset$-definable family $(B_z)_{z\in Z}$ of balls $B_z\subseteq K^n$, and a $\emptyset$-definable family $(X_{z})_{z\in Z}$ of subsets $X_z \subseteq K^n$, all parametrized by a $\emptyset$-definable
set $Z$.
Then the map
\[
Z \to \bigcup_d G^d_n(k), z \mapsto \rtsp_{B_z}(X_{z})
\]
is $\emptyset$-definable. Moreover, a formula defining that map in $K$ also works in any other spherically complete $K' \equiv K$.
\end{thm}

Another important ingredient we will need is that one can control precisely on which balls a definable set is not riso-trivial.

\begin{thm}\label{t.S0}
Given a definable set $X \subseteq K^n$,
there exist finitely many disjoint sets $B_1, \dots, B_\ell \subseteq K^n$ each of which is either a singleton or a ball such that for every ball $B \subseteq K^n$, $X$ is not riso-trivial on $B$ if and only if $B_i \subseteq B$ for some $1 \le i \le \ell$.
\end{thm}

\begin{proof}
By \cite[Lemma~3.4.3]{iBW.canStrat}, the ``rigid core'' $\operatorname{Crig}_X$ is a finite set of pairwise disjoint singletons and balls. By the definition of $\operatorname{Crig}_X$ (just before that lemma), we have $\operatorname{Crig}_X = \{B_1, \dots, B_\ell\}$.
\end{proof}

Note that if $B$ is a proper sub-ball of some of the $B_i$, then necessarily $X$ is $1$-riso-trivial on $B$, since then, $B$ contains no $B_j$.

\begin{rem}\label{r.S0.uniq}
The sets $B_1, \dots B_\ell$ are uniquely determined by $X$, up to permutation, namely
these are
\begin{enumerate}
 \item 
those balls $B$ for which $X$ is not riso-trivial on $B$ but $1$-riso-trivial on any proper subball of $B$; and
\item those singletons $\{b\}$ such that $X$ is not riso-trivial on every ball containing $b$.
\end{enumerate}
\end{rem}

\begin{rem}\label{r.S0.bound}
By Theorem~\ref{t.defble}, the characterization of the $B_i$ given in Remark~\ref{r.S0.uniq} is definable: If $(X_{z})_{z\in Z}$ is a $\emptyset$-definable familiy of sets $X_z \subseteq K^n$ and $C_z$ is the union of the points and balls $B_{z,1}, \dots, B_{z,\ell_z}$ corresponding to $X_z$, then $(C_z)_z$ is also a $\emptyset$-definable family. In addition, for any such family, there exists $N \in \NN$ bounding the cardinalities $\ell_z$ for all $z \in Z$.
(This follows by applying \cite[Lemma~3.4.3]{iBW.canStrat} in a language with a constant symbol for $z$ added.)
\end{rem}

\begin{ex}
Set $X := \OK \times \{0\} \subset K^2$. Then clearly
$X$ is $1$-riso-trivial on $\OK \times \OK$ but not on $B := t^{-1}\OK \times t^{-1}\OK$ (recall that $t \in \OK$ is a generator of the maximal ideal of $K$). In this way, one obtains that the finitely many sets $B_i$'s provided for $X$ by Theorem~\ref{t.S0} consist only of $B$.
\end{ex}

\begin{ex}
Let $X \subset K^2$ be the graph of the function $x \mapsto x^2$. Then one can verify that $X$ is not riso-trivial on $\OK \times \OK$ but $1$-riso-trivial on each proper subball of $\OK \times \OK$. This shows that $B := \OK \times \OK$ is among the balls provided by Theorem~\ref{t.S0}. One then also checks that $X$ is $1$-riso-trivial on every ball disjoint from $B$, so $B$ is the only ball given by Theorem~\ref{t.S0} for $X$.
\end{ex}

\begin{figure}[hbt]
    \centering
    \includegraphics[scale=0.65]{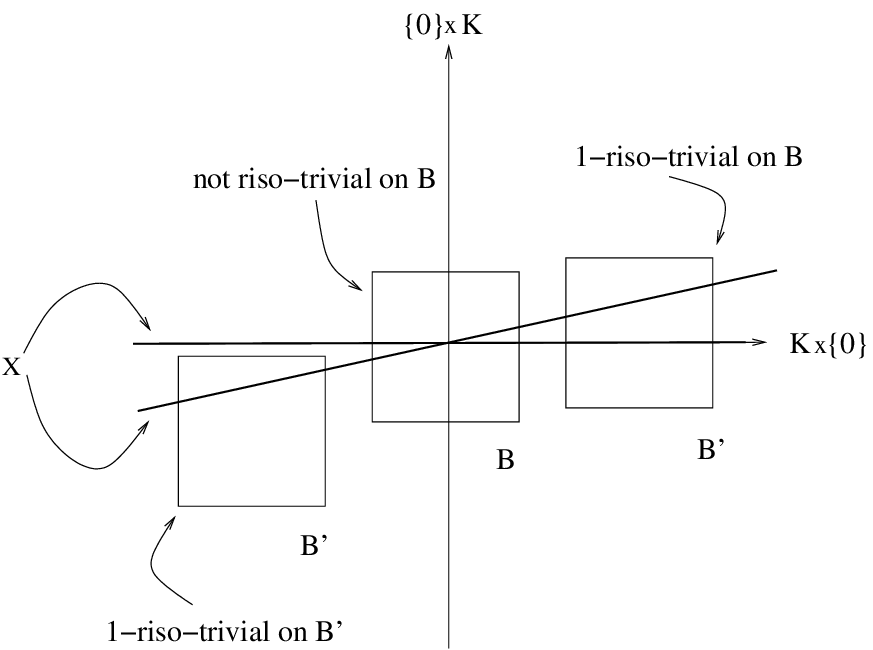}
    \caption{Applying Theorem~\ref{t.S0} in Example~\ref{e.two-lines}}
    \label{f.two-lines}
\end{figure}

\begin{ex}\label{e.two-lines}
Let $X \subset K^2$ be the union of the lines $K \times \{0\}$ and $\{(x,tx);  x \in K\}$ (see Figure~\ref{f.two-lines}). 
Then $X$ is not riso-trivial on every ball $B$ containing the origin. If $B' \subset K^2$ is a ball not containing the origin, then $X$ is clearly $1$-riso-trivial on $B'$ if $B'$ meets only one of the two lines, but even if $B'$ meets both lines, $X$ is $1$-riso-trivial on $B'$, since the two lines can be made parallel using a risometry. Overall, we obtain that Theorem~\ref{t.S0} yields only the set $\{(0,0)\}$.

\end{ex}

\section{Motivic integrals}\label{s.Motivic integrals}

In this section, we give a brief review on motivic integration, in order to fix notation and conventions 
(for all the details concerning this theory see \cite{CL.mot} or \cite{CL.mot0p}). In addition, we introduce the notion of limits of constructible motivic functions and prove some basic properties, adapting ideas from \cite[Section 3.1]{For.motCones}.

In \cite{CL.mot0p} motivic integration is introduced under axiomatic assumptions stated as \cite[Definition~3.9]{CL.mot0p}. By \cite[Theorem~5.8.2]{iCR.hmin}, those assumptions are satisfied for $\cL$ and $K$ as in Definition~\ref{d.cL}.

\begin{notn}  
\label{not:motivic integration} 
For a valued field $K$ with residual field $k$ and value group $\ZZ$ as in Section \ref{s.Language}, let us first recall the definition given in \cite{CL.mot} of the ring of constructible motivic functions on some $\emptyset$-definable set $X\subseteq K^n\times k^m\times \ZZ^r$. 
For this consider $\mathrm{RDef}_X$
the set of $\emptyset$-definable subsets of $X\times k^p$, for $p\in \NN$, equipped with their natural projection on $X$. 
One denotes by $\mathcal Q_+(X)$ the Grothendieck semiring of $\emptyset$-definable isomorphism classes $[Y]$  (more accurately classes of projections from $Y$ to $X$), $Y\in \mathrm{RDef}_X$, with the following relations for any $\emptyset$-definable sets 
$Y,Z\in \mathrm{RDef}_X$
\begin{itemize}
\item[-] $[\emptyset]=0$,
\item[-] $[Y\cup Z]+[Y\cap Z]=[Y]+[Z]$, 
\item[-] $[Y\times Z]=[Y]\cdot[Z]$.
\end{itemize}
We denote by the formal symbol $\LL$, the same symbol as the one used in Notation \ref{n.norm} to define the norm, the class $[X\times k]\in \mathcal Q_+(X)$ (actually the class of the projection of $X\times k$ on $X$), and by $\mathcal Q(X)$ the associated abelian Grothendieck ring. 

Still with the same symbol $\LL$  we consider on the other hand the ring 
$$
A:=\ZZ[\LL, \LL^{-1}, \left(\frac{1}{1-\LL^i}\right)_{i>0}],
$$ 
and $A_+$ its semigroup of nonnegative elements, where positivity is defined by evaluations 
$\LL\mapsto r$ at $r$, for any real number $r>1$. 
One defines the ring ${\mathcal P}(X)$ of Presburger functions, as a subring of the ring of functions 
$X\to A$, as generated by constant functions, $\emptyset$-definable functions $X\to \ZZ$ and functions $\LL^\beta$, with $\beta: X\to \ZZ$ $\emptyset$-definable. 
The semiring of functions of  ${\mathcal P}(X)$ with nonnegative values (that is with values in $A_+$) is denoted ${\mathcal P}_+(X)$. 

We denote by ${\mathcal P}_+^0(X)$ (resp. ${\mathcal P}^0(X)$) the sub-semiring (resp. sub-ring) of ${\mathcal P}_+(X)$ (resp. ${\mathcal P}$(X))
generated by all characteristic functions $\mathbf{1}_Y$, with $Y\subseteq X$, and the constant function $\LL-1$. We can send ${\mathcal P}_+^0(X)$ (resp. ${\mathcal P}(X)$)
in $\mathcal Q_+(X)$ (resp. $\mathcal Q(X)$) by 
$\mathbf{1}_Y \mapsto [Y] $ and 
$\LL-1\mapsto \LL-1$, to form the semiring $\mathcal{C}_+(X)$ of nonnegative constructible motivic functions on $X$ in the following way 
$$
\mathcal{C}_+(X):=\mathcal Q_+(X)\otimes_{{\mathcal P}_+^0(X)} {\mathcal P}_+(X),
$$
and likewise the ring $\mathcal{C}(X)$ of constructible motivic functions on $X$ in the following way 
$$
\mathcal{C}(X):=\mathcal Q(X)\otimes_{{\mathcal P}^0(X)} {\mathcal P}(X).
$$

\begin{rem}\label{r.positivity in P+}
One can define (see \cite[Sections 4.2, 12.2]{CL.mot}) a partial order on ${\mathcal P}(X)$, denoted $\le$, defined by   $$\forall \varphi, \psi \in {\mathcal P}(X), \ \varphi\le \psi \Longleftrightarrow \psi-\varphi \in {\mathcal P}_+(X).$$
Similarly, one can define a partial preorder  on $\mathcal{C}(X)$, still denoted $\leq$, by replacing $\mathcal{P}(X)$ by $\mathcal{C}(X)$ and $\mathcal{P}_+(X)$ by $\mathcal{C}_+(X)$
in the above definition. 
In both situation, when $\varphi\le \psi$, one says that $\varphi$ is bounded by $\psi$. 
\end{rem}

With this notation, in \cite[Theorem 10.1.1]{CL.mot} the functor of motivic integration 
is defined as a functor from the category of $\emptyset$-definable sets (equipped with their projection on $X$) to the category of abelian semigroups, sending a $\emptyset$-definable set $Y$ that projects on $X$ to the semigroup $\mathrm{I}_X\mathcal{C}_+(Y)\subseteq \mathcal{C}_+(Y)$ of nonnegative integrable functions on $Y$ over $X$. To any $\emptyset$-definable map $f:Y\to Z$ between $\emptyset$-definable sets of $K^n\times k^m\times \ZZ^r$, it associates a morphism
$$
f_!: \mathrm{I}_X \mathcal{C}_+(Y) \to \mathrm{I}_X \mathcal{C}_+(Z),
$$
corresponding to integration in the fibres of $f$ twisted by jacobian, in the sense that for $f:Y\to Z$ and $g:Z \to W$, and $\varphi\in \mathrm{I}_X\mathcal{C}_+(Y)$,
$(g\circ f)_!(\varphi)=g_!(f_!(\varphi))$.  

\begin{rem}[{\cite[Theorem 10.1.1 A0 (c)]{CL.mot}}]\label{r.Fubini}
One has  $\varphi\in \mathrm{I}_X\mathcal{C}_+(Y)$ if and only if 
$\varphi\in \mathrm{I}_Z\mathcal{C}_+(Y)$ and $f_{!}(\varphi)\in \mathrm{I}_X\mathcal{C}_+(Z)$.
\end{rem}

In case $X=Z=\{pt\}$, we have
$$ \mathrm{I}_{X} \mathcal{C}_+(Z)=\mathrm{I}_{\{pt\}} \mathcal{C}_+(\{pt\})= \mathcal{C}_+(\{pt\})=\mathcal Q_+(\{pt\})\otimes_{\NN[\LL-1]} A_+,$$ 
where $\mathcal Q_+(\{pt\})$ is the Grothendieck semiring of $\emptyset$-definable subsets of $k^m$, for $m\in \NN$.
In this situation, $f_!(\varphi)$, denoted $\mu_n(\varphi)$ or $\mu(\varphi)$ when the context is clear, is the motivic integral of $\varphi$ on $Y$, $\varphi$ is said integrable on $Y$ and we denote $\mathrm{I}\mathcal{C}_+(Y)$ instead of $\mathrm{I}_{\{pt\}}\mathcal{C}_+(Y)$. In case moreover $\varphi=\mathbf{1}_Y$, $\mu_n(\varphi)$ is simply denoted $\mu_n(Y)$, or $\mu(Y)$, and is called the motivic volume of $Y$.

\begin{rem}\label{r.normalisation}
A normalisation of the volume is made in order to obtain uniqueness of the functor of motivic integration, leading in particular to $\mu_1(\OK)=1$ and 
$\mu_1(\mathcal{M}_K)=\LL^{-1}$ (see \cite[Theorem 10.1.1, A7]{CL.mot}). 
\end{rem}

We mention 
without details that one can build motivic measures in all dimensions, that is to say measures $\mu_d$ adapted to functions with support a $\emptyset$-definable set of $K^n$ of dimension $d$, with $d\in \NN$ (see \cite[Section 6]{CL.mot}, \cite[Section 11]{CL.mot0p}).


\begin{rem}\label{r.subset and inequality}
Let  $X,Y$ be $\emptyset$-definable subsets of $K^n\times k^m \times \ZZ^r$, of dimension $d$, say, such that $X\subset Y$. Then, according to Remark \ref{r.positivity in P+}, $\mathbf{1}_X$ is bounded by $\mathbf{1}_Y$ in $\mathcal{P}_+(Y)$ (or in $\mathcal{C}_+(Y)$), since $\mathbf{1}_Y-\mathbf{1}_X=\mathbf{1}_{Y\setminus X}\in \mathcal{P}_+(Y)$,
and thus integrating, by linearity of the integral and preservation of nonnegativity, one gets $\mu_d(X)\le \mu_d(Y)$.
\end{rem}

In what follows we mainly consider the family version of this construction (see \cite[Section 14]{CL.mot}, or \cite{CL.mot0p} for a specific construction of integration in fibres in the more general context of residue fields of all characteristics), that is,  for some integers $r,s,t$, we consider $\Lambda\subset K^r\times k^s\times \ZZ^t$, a $\emptyset$-definable set, and $\emptyset$-definable families
$Y=(Y_\lambda)_{\lambda\in \Lambda}$, $Z=(Z_\lambda)_{\lambda\in \Lambda}$ compatible with the projection on $X$, in the sense that
$Y$ and $Z$ are given by maps $Y\to X\to \Lambda$ and $Z\to X\to \Lambda$.

\begin{notn}
For $n,m,r\in \NN$, $Y,Z\subseteq K^n\times k^m\times \ZZ^r$ $\emptyset$-definable sets, $f:Y\to Z$ a $\emptyset$-definable map and $\varphi\in \mathcal{C}(Z)$ we denote
by $f^*(\varphi)\in \mathcal{C}(Y)$ the pull-back of $\varphi$ along $f$, which should be understood as the composition $\varphi\circ f$. 
In case $f$, given by a commutative diagram
$$
  \xymatrixcolsep{1pc}\xymatrix{
    Y   \ar[dr]  \ar[rr]^{f}  & &
    Z \ar[ld]  &\\
 & X \ar[d]  & \\
 & \Lambda &
  } 
 $$
 for $\lambda\in \Lambda$, we denote by $f_\lambda:Y_\lambda\to Z_\lambda$ the restriction of $f$ to the fibres of $Y$ and $Z$ above $\lambda$. We also denote, with a slight abuse of notation, $i_\lambda $  the inclusion of $Y_\lambda$ (resp. of $Z_\lambda$) in $Y$ (resp. in $Z$).
 \end{notn}
 
 One can then define (\cite[Theorem 14.1.1]{CL.mot})
a morphism
$$
f_{!\Lambda}: \mathrm{I}_X \mathcal{C}_+(Y\to \Lambda) \to \mathrm{I}_X \mathcal{C}_+(Z\to \Lambda),
$$
that associates to any $\varphi$ in the  semigroup $\mathrm{I}_X \mathcal{C}_+(Y\to \Lambda)$ of integrable nonnegative constructible motivic functions, graded with respect to relative dimension (relative to the map $Y\to X$), the element $f_{!\Lambda}(\varphi)$ of $\mathrm{I}_X \mathcal{C}_+(Z\to \Lambda)$.

\begin{rem}\label{r.notation famille}
By \cite[Proposition 14.2.1]{CL.mot} one has for any 
$\lambda \in \Lambda$,  $$i^*_\lambda(f_{!\Lambda}(\varphi))=f_{\lambda!}(i^*_\lambda(\varphi)),$$
This is why we will often denote $f_{!\Lambda}(\varphi)$ by 
$(f_{\lambda!}(\varphi_\lambda))_{\lambda\in \Lambda}$, that is to say that we will consider the constructible motivic function $f_{!\Lambda}(\varphi)$ as a family of measures over $\Lambda$. 
\end{rem}

In the particular case $f:Y\to \Lambda$, that is $Z=X=\Lambda$ and the map $f$ is the map defining the family $Y$, we have  $f_\lambda : Y_\lambda\to \{\lambda\}$, and consequently
$f_{!\Lambda}(\varphi)=(\mu(\varphi_\lambda))_{\lambda\in \Lambda}\in \mathcal{C}_+(\Lambda)$. 
In this situation $f_{!\Lambda}(\varphi)$ is the integration of 
$\varphi$ in the fibres of $Y\to \Lambda$, and
we denote $f_{!\Lambda}$ by $\mu_\Lambda$, or in a more intuitive way by
$\left( \int_{y\in Y_\lambda} \varphi_\lambda(y) \ dy \right)_{\lambda\in \Lambda}$. When furthermore $\varphi=\mathbf{1}_Y$, $\mu_\Lambda(\varphi)=(\mu(Y_\lambda))_{\lambda\in \Lambda}$.

\end{notn}

\begin{rem}\label{rem.direct product} 
In the situation where $\pi: X\times Y\to X$ is the standard projection, then $\pi_!=\pi_{!X}$, and thus 
$\pi_!(\varphi)\in \mathrm{I}_X\mathcal{C}_+(X)=\mathcal{C}_+(X)$ is the family $\left(\mu(\varphi_x)\right)_{x\in X}=\left(\int_{Y} \varphi(x,y) \ dy\right)_{x\in X}$  of motivic integral of $\varphi$ in the fibres $\{x\}\times Y$, for $x\in X$ (see \cite[Remarque 5.3.4.6.1]{RA.hdr}).
\end{rem}

We now introduce a notion of convergence of nonnegative construc\-ti\-ble motivic functions and measures. A variant of this has already been considered in \cite[Section 3.1]{For.motCones}, the difference here being that we do not take mean values in the notion of convergence.

One can uniquely extend the degree (in $\LL$) on $\ZZ[\LL]$ to a degree $\deg$ on $A_+$ (resp. on $A$). We then define a topology on $A_+$ such that the convergence to $0$ of a sequence $(a_n)_{n\in \NN}$ with elements in $A_+$ is the convergence to $-\infty$ of the sequence $(\deg(a_n))_{n\in \NN}$. This yields a notion of convergence of functions in $\mathcal{P}_+(\NN)$, and also in a variant with additional parameters.

\begin{defn}
Let $X,Y\subset K^n\times k^m\times \ZZ^r$ be $\emptyset$-definable, and let $\psi \in \mathcal{P}_+(X \times Y)$ and $\psi' \in \mathcal{P}_+(X)$ be nonnegative Presburger function.
\begin{enumerate}
    \item 
In the case $Y = \NN$,
we say that $\psi(\cdot,r)$ \emph{ converges to $\psi'$, for $r \to +\infty$}, if for each $x \in X$, we have $\lim_{r \to +\infty}\psi(x,r) = \psi'(x)$, in the sense of convergence in $A_+$ described above.
  \item In the case $Y = \NN$, we say that $\psi(\cdot,r)$ is \emph{increasing for $r \to +\infty$} if for each $(x,r)\in X \times \NN$, we have $\psi(x,r+1) - \psi(x,r) \in A_+$.
  \item We say that $\psi(x, y)$ is \emph{bounded by $\psi'$ (uniformly in $y$)}, if for each $(x,y)\in X \times Y$, we have $\psi'(x) - \psi(x,y) \in A_+$ (see also Remark \ref{r.positivity in P+}).
\end{enumerate}
\end{defn}

We now extend those notions to constructible motivic functions.
By orthogonality of the residue field and the value group, any function $\psi\in \mathcal{C}_+(X \times \NN)$ may be written as a finite sum
\begin{equation}\label{eq.decomposition}
\psi =\sum_{i=1}^\ell c_i\otimes \psi_i,
\end{equation} 
where $c_i\in \mathcal Q_+(X)$ and 
$\psi_i\in \mathcal{P}_+(X \times \NN)$ (see \cite[Section 5.3]{CL.mot} or \cite[Proposition 7.5]{CL.mot0p} for details).

\begin{defn}\label{d.convergence Mot}
Let $X\subset K^n\times k^m\times \ZZ^r$ be $\emptyset$-definable and let $\psi \in \mathcal{C}_+(X \times \NN)$ and $\psi' \in \mathcal{C}_+(X)$ be nonnegative constructible motivic functions.
We say that $\psi(\cdot,r)$ \emph{converges to $\psi'$ (resp.\ is increasing, resp.\ is bounded by $\psi'$)} if
there exist $c_i$ and $\psi_i$ as in \eqref{eq.decomposition} and $\psi'_i \in  \mathcal{P_+}(X)$ such that $\psi' = \sum_i c_i\otimes\psi'_i$ and such that for each $i$, $\psi_i$ converges to $\psi'_i$ (resp.\ is increasing resp.\ is bounded by $\psi'_i$). For the notion of bounded, we also allow arbitrary definable set $Y\in K^n\times k^m \times \ZZ^r$ instead of $\NN$.
\end{defn}

\begin{rem}
The notions of bounded and increasing can equivalently be defined directly using $\mathcal{C}_+$, that is without passing through $\mathcal P_+$ (see also Remark \ref{r.positivity in P+}): $\psi \in \mathcal{C}_+(X \times Y)$ is bounded by $\psi' \in \mathcal{C}_+(X)$ if there exists $\phi \in  \mathcal{C}_+(X \times Y)$ such that $\psi + \phi = \psi'$, and $\psi \in \mathcal{C}_+(X \times \NN)$ is increasing if there exists $\phi \in  \mathcal{C}_+(X \times \NN)$ such that $\psi + \phi = \psi'$, where $\psi'(x,r) = \psi(x,r+1)$.
\end{rem}

\begin{notn}\label{n.le}
For $n,m,r\in \NN$, $X\subseteq K^n\times k^m\times \ZZ^r$ an $\emptyset$-definable set and $\phi,\psi \in \mathcal C_+(X)$ we write $\phi \le \psi$ if $\phi$ is bounded by $\psi$ in the sense of Definition \ref{d.convergence Mot}.
\end{notn}

\begin{rem}\label{r.fonct}
The limit $\psi'$, if it exists, is well-defined, that is does not depend on the choice of the representation of $\psi$ in the form \eqref{eq.decomposition}. To see this denote by $\mathcal P^{\mathrm{lim}}_+(X \times \NN) \subset \mathcal P_+(X \times \NN)$ and $\mathcal{C}^{\mathrm{lim}}_+(X \times \NN) \subset \mathcal{C}_+(X \times \NN)$ the sub-semirings of converging nonnegative Presburger and constructible motivic functions, where the convergence is defined in Definition \ref{d.convergence Mot}. By definition of $\mathcal{C}^{\mathrm{lim}}_+(X \times \NN)$, we have
$$\mathcal{C}^{\mathrm{lim}}_+(X \times \NN)=
\mathcal Q_+(X)\otimes_{\mathcal P_+^0(X)} \mathcal P^{\mathrm{lim}}_+(X \times \NN),$$
where we use again the same isomorphism as in \eqref{eq.decomposition}, provided by 
\cite[Section 5.3]{CL.mot} or \cite[Proposition 7.5]{CL.mot0p}, to omit the $\NN$ from the left side of the tensor product and from the $\mathcal P_+^0$.

The map $h\colon \mathcal{P}^{\mathrm{lim}}_+(X \times \NN) \to \mathcal{P}_+(X)$ sending a function to its limit induces the identity on
$\mathcal P_+^0(X)$, so we can apply the functor
$S \mapsto \mathcal Q_+(X)\otimes_{\mathcal P_+^0(X)} S$
to that map $h$ to define the limit on $\mathcal C^{\mathrm{lim}}_+(X \times \NN)$ in a way which is independent of how we write a function $\psi \in \mathcal C^{\mathrm{lim}}_+(X \times \NN)$ as in \eqref{eq.decomposition}.

\end{rem}

The following Lemma \ref{l.squeeze} and Remarks \ref{r.monoton convergence}, \ref{r.tub}
give the properties of this notion of convergence that we need. The proof of Lemma \ref{l.squeeze} needs a detailed understanding of the tensor product $\mathcal Q_+(X) \otimes_{\mathcal{P}^0_+(X)} \mathcal{P}_+(X \times \NN)$. We give those details in an appendix (Section~\ref{s.tensor}); more precisely, we use Lemma~\ref{l.trivial-writing} from the appendix in the proof of Lemma \ref{l.squeeze}.


\begin{lem}\label{l.squeeze}
Suppose that we have $\psi_1, \psi_2\in \mathcal{C}_+(X \times \NN)$ such that $\psi_1(\cdot, r) + \psi_2(\cdot, r)$ tends to $0 \in \mathcal{C}_+(X)$ (for $r \to +\infty$). Then 
$\psi_1(\cdot, r)$ and $\psi_2(\cdot, r)$ also tend to $0$ for $r \to +\infty$.
\end{lem}

Note that most likely, this lemma would be false if, instead of working with $\mathcal{C}_+(X \times \NN)$ and $\mathcal{C}_+(X)$, we would work with the image of these semirings in the rings $\mathcal{C}(X \times \NN)$ and in $\mathcal{C}(X)$, respectively. Indeed, it seems plausible that for sufficiently bad residue fields, one can find some $\phi \in \mathcal{C}(X) \setminus \{0\}$ satisfying $2\phi = 0$ and which lies in the image of $\mathcal{C}_+(X)$. In other words, we would not even have that $\psi_1+\psi_2 = 0$ implies $\psi_1 = \psi_2 = 0$.

\begin{proof}[Proof of Lemma~\ref{l.squeeze}]
We first collect some facts that will serve as prerequisites to apply Lemma~\ref{l.trivial-writing}. All those facts follow directly from the definitions of the corresponding semirings. Given $\psi \in P_+(X \times \NN)$, we write $\psi \to 0$ to say that $\psi(\cdot, r)$ converges to $0$ when $r$ goes to $+\infty$ (in the sense of Definition \ref{d.convergence Mot}).
\begin{enumerate}
    \item For every $\emptyset$-definable set $Y$ and every $\phi, \phi' \in \mathcal Q_+(Y)$, we have the following two implications:\\(1.a) $\phi+\phi' = 0$ implies $\phi = \phi' = 0$,
    \\(1.b) $\phi\cdot \phi' = 0$ implies $\phi = 0$ or $\phi' = 0$.
    \item For every $\emptyset$-definable set $Y$ and every $\phi, \phi' \in \mathcal{P}_+(Y)$, we have the following two implications:\\(2.a) $\phi+\phi' = 0$ implies $\phi = \phi' = 0$,\\(2.b) $\phi\cdot \phi' = 0$ implies $\phi = 0$ or $\phi' = 0$.
    \item For every $\phi_1,\phi_2 \in \mathcal P_+(X \times \NN)$ and every $\phi_0 \in \mathcal P_+^0(X)$, we have the following two implications:\\
    (3.a) $\phi_1 + \phi_2 \to 0$ implies
    $\phi_1 \to 0$ and $\phi_2 \to 0$,\\ (3.b)
    $\phi_0 \cdot \phi_1 \to 0$ implies $\phi_1 \to 0.$
\end{enumerate}\label{eq.++}
Suppose now that $\psi_1, \psi_2$ are given as in the lemma. Write both of them as sums as in \eqref{eq.decomposition}:
\begin{equation}
\psi_1 =\sum_{i=1}^{\ell_1} c_{1,i}\otimes \psi_{1,i},
\qquad
\psi_2 =\sum_{i=1}^{\ell_2} c_{2,i}\otimes \psi_{2,i},
\end{equation} 
for some $c_{j,i} \in \mathcal Q_+(X)$ and 
$\psi_{j,i}\in \mathcal{P}_+(X \times \NN)$.
By assumption, $\psi_3 := \psi_1 + \psi_2$ converges to $0$, i.e. it can also be written as such a sum
\begin{equation}\label{eq.psi3}
\psi_3 =\sum_{i=1}^{\ell_3} c_{3,i}\otimes \psi_{3,i},
\end{equation} 
where each $\psi_{3,i}$ converges, say to some $\phi_{3,i} \in \mathcal{P}_+(X)$,
and such that the sum
\begin{equation}\label{s.phi0}
\sum_{i=1}^{\ell_3} c_{3,i}\otimes \phi_{3,i} \end{equation} 
is equal to $0$ in $\mathcal{C}_+(X) = \mathcal Q_+(X) \otimes_{\mathcal P_+^0(X)} \mathcal P_+(X)$. We may without loss assume $c_{3,i} \ne 0$ for each $i$.

Apply Lemma~\ref{l.trivial-writing} to
$S := \mathcal P_+^0(X)$,
$U_1 := \{0\} \subset \mathcal Q_+(X)$ and $U_2 := \{0\} \subset \mathcal P_+(X)$. Note that the assumptions of Lemma \ref{l.trivial-writing} concerning $S$, $U_1$ and $U_2$ are satisfied by (2), (1) and (2), respectively.
Since the sum \eqref{s.phi0} lies in the sub-semimodule $U \subset \mathcal Q_+(X) \otimes_{\mathcal P^0_+(X)} \mathcal P_+(X)$ defined in Lemma \ref{l.trivial-writing}, we obtain $\phi_{3,i} = 0$ for each $i$ (since we assumed $c_{3,i} \ne 0$). In other words, for each $i$, we have $\psi_{3,i} \to 0$.

Now apply Lemma~\ref{l.trivial-writing} once more, this time to
$S := \mathcal P_+^0(X)$,
$U_1 := \{0\} \subset \mathcal Q_+(X)$ and $U_2 := \{\phi \in \mathcal P_+(X \times \NN); \phi \to 0\}$. The assumptions of Lemma \ref{l.trivial-writing} concerning $S$, $U_1$ and $U_2$ are satisfied by (2), (1) and (3). Since in \eqref{eq.psi3}, $\psi_{3,i} \to 0$ for each $i$, $\psi_3$ lies in the sub-semimodule $U \subset \mathcal Q_+(X) \otimes_{\mathcal P^0_+(X)} \mathcal P_+(X \times \NN)$ defined in Lemma \ref{l.trivial-writing}, so the lemma implies that $\psi_1$ and $\psi_2$ lie in $U$. By definition of $U$, elements of $U$ converge to $0$, so we are done.
\end{proof}

\begin{rem}[Monotone convergence]\label{r.monoton convergence}
A $\emptyset$-definable subset $X$ of $K^n\times \NN$ with projection $\pi: X \mapsto  \NN$, and fibres $\pi^{-1}(p)=X_p$, such that $(X_p)_{p\in \NN}$ is uniformly bounded with respect to $p\in \NN$, and which is an increasing sequence with respect to the inclusion relation, gives rise to 
the bounded increasing function $\varphi=\pi_{!\NN}(\mathbf{1}_{X})=(\mu(X_p))_{p\in \NN}\in \mathcal{C}_+(\NN)$.
It follows that this function converges as $p$ goes to $+\infty$, and by definition of the motivic integral, that its limit is $\mu\left( \bigcup_{i=0}^{+\infty} X_i\right)$ (see \cite[Lemma 5.1]{For.motCones} for more details). \end{rem}

\begin{rem}\label{r.tub}
Let $X \subset K^n$ be a bounded $\emptyset$-definable set of dimension $d$ and for $r \in \ZZ$, let
$T_r(X)=\{x\in K^n;\exists y\in X, \val(y-x)\ge r\}$ be the corresponding tubular neighbourhood. Then, for $r \ge 0$,
$\mu_n(T_r(X))$ is bounded by $(\LL^{-r})^{n-d}\cdot C$ for some constant $C \in \mathcal{C}_+(\{pt\})$ depending only on $X$ (i.e. $(\LL^{-r})^{n-d}\cdot C - \mu_n(T_r(X)) \in \mathcal{C}_+(\NN)$). Indeed, by cell decomposition with $1$-Lipschitz centres, $X$ can be partitioned into a finite $\emptyset$-definable family of sets $X_\xi$ each of which is, up to a permutation of coordinates, the graph of a $1$-Lipschitz function $f_\xi$ on a subset $\bar X_\xi$ of $K^d$ (for more details about this kind of cell decomposition, see Step~1 of the proof of Theorem~\ref{t.Crofton}). This implies
that $\mu_n(T_r(X_\xi))$ is bounded by
$\mu_d(T_0(\bar X_\xi))\cdot (\LL^{-r})^{n-d}$. Finally note that, 
applying Remark \ref{r.subset and inequality},
the volume of the tubular neighbourhood of $X$ is bounded by the sum of the volume of the tubular neighbourhoods of the $X_\xi$'s. 
\end{rem}
\section{A partial preorder for constructible motivic functions}\label{s.preorder}

In this section we define and study a partial preorder $\preccurlyeq$ on the set of integrable constructible motivic functions. We show that $\preccurlyeq$ has the required properties with respect to our proof strategy of the nonarchimedean version of inequality \eqref{eq.EntPol} relating the motivic Vitushkin invariants and the motivic metric entropy. In particular $\preccurlyeq$ is compatible with integration (Lemma \ref{l.<int}). We first start defining the set 
$\mathrm{I}\mathcal{C}_{\succcurlyeq0}(X)$ of constructible motivic functions that are nonnegative with respect to our new relation, and observe that nonnegative 
constructible motivic functions in the sense of Section \ref{s.Motivic integrals} (that is to say when nonnegativity is defined through the ring $A_+$) are also nonnegative regarding our definition of $\preccurlyeq$. 
Note that  with respect to our preorder $\preccurlyeq$, there are strictly more nonnegative functions than with respect to the preorder $\le$ defined by $\mathrm{I}\mathcal{C}_{+}(X)$ (see Example~\ref{e.more}).
\begin{defn}\label{d.<}
Let $m,n,r\in \NN$ and $X\subseteq K^n\times k^m\times \ZZ^r$ be a $\emptyset$-definable set. We denote by $\mathrm{I}\mathcal{C}_{\succcurlyeq0}(X) $ the subset of $\mathrm{I}\mathcal{C}(X)$ consisting of the constructible motivic functions $F\in \mathrm{I}\mathcal{C}(X)$ for which
there exist $\emptyset$-definable sets $Y, Z \subseteq k^{m'}$ for some $m'\in \NN$, a $\emptyset$-definable surjection $f\colon Z \twoheadrightarrow Y$ with finite fibres and a constructible motivic function $\varphi \in \mathrm{I}\mathcal{C}_+(Y \times X)$ such that the following equality holds in $\mathrm{I}\mathcal{C}(X)$
 \begin{equation}\label{eq.d.<}
 F = \pi_{X!}f_X^*\varphi - \pi_{X!}\varphi,
 \end{equation}
 where $f_X: Z \times X \twoheadrightarrow Y \times X, (z,x) \mapsto (f(z), x)$, and $\pi_X$ is the standard projection onto $X$.
$$
  \xymatrixcolsep{2pc}\xymatrix{
    Z\times X  \ar@{->>}[dr]_{\pi_X} \ar@{->>}[rr]^{f_X}  &   & 
    Y\times X \ar@{->>}[ld]^{\pi_X} 
    \\
 & X  \ar@{->>}[d] & \\
 & \{ pt \}  & 
  } 
 $$
 
\end{defn}

\begin{rem}
Note that in \eqref{eq.d.<}, $\varphi$ and $f^*_X\varphi$ are indeed integrable along $\pi_X$:
For $\varphi$, this follows from the assumption $\varphi \in \mathrm{I}\mathcal{C}_+(Y \times X)$ (which implies $\varphi \in \mathrm{I}_{X}\mathcal{C}_+(Y \times X)$ by Remark \ref{r.Fubini}.
For $f^*_X\varphi$, note that since $f$ has finite fibres, $\varphi \in \mathrm{I}_X\mathcal{C}_+(Y \times X)$ implies $f^*_X\varphi \in  \mathrm{I}_X\mathcal{C}_+(Z \times X)$.
\end{rem}

In view of Remark \ref{rem.direct product}, we have the more intuitive writing for $F\in \mathrm{I}\mathcal{C}_{\succcurlyeq 0}(X)$. 
 $$
 F(x)=\mu_X(f_X^*\varphi_x)-\mu_X(\varphi_x)=
 \int_{Z} \varphi\circ f_X(z,x) \ dz - \int_{Y} \varphi(y,x) \ dy.
$$

One could imagine various variants of Definition~\ref{d.<}; we will see in Lemma~\ref{l.<surj} and Remark~\ref{r.<sorts} below that some variants lead to an equivalent definition.

\begin{lem}\label{l.<C+}
For $m,n,r\in \NN$ and $X\subseteq K^n\times k^m\times \ZZ^r$ a $\emptyset$-definable set one has that the image of $\mathrm{I}\mathcal{C}_+(X)$ in $\mathrm{I}\mathcal{C}(X)$ is contained in $\mathrm{I}\mathcal{C}_{\succcurlyeq0}(X)$.
\end{lem}
\begin{proof}
Let $F \in \mathrm{I}\mathcal{C}_+(X)$ be given. Set $Y := \{0\}\subseteq k$, $Z = \{0, 1\}\subseteq k$, $f(0) = f(1) = 0$ and  $\varphi(0,\cdot) = F$.
Then $\pi_{X!}(f_X^*\varphi) - \pi_{X!}(\varphi) = 2F - F = F$.
\end{proof}

\begin{ex}\label{e.more}
In case $k=\CC$ for instance, one can give an example of a constructible motivic 
function in $\mathrm{I}\mathcal{C}_{\succcurlyeq0}(X)$ but not in the image of $\mathrm{I}\mathcal{C}_+(X)$ in $\mathrm{I}\mathcal{C}(X)$. Indeed, let us consider the constructible motivic function $\varphi=[Z]-1$ for $Z = \{\pm \sqrt2\} \subset \CC$. In view of Lemma \ref{l.<surj}, 
since there is a $\emptyset$-definable surjection from the $\emptyset$-definable set $Z$ to the point, we have $\varphi\in \mathrm{I}\mathcal{C}_{\succcurlyeq0}(X)$. 

Let us now show that $\varphi$ does not lie in the image of $\mathrm{I}\mathcal{C}_+(X)$. This amounts to showing that there is no $\emptyset$-definable set $Y \subset \CC^m$ (for some $m$) and injective $\emptyset$-definable map $f\colon \{pt\} \sqcup Y \to Z \sqcup Y$. Suppose $Y$ and $f$ would exist. Without loss, they are defined by quantifier free formulas, so that we also obtain a $\emptyset$-definable injection between the restrictions of those sets to any subfield $F \subset \CC$. We pick $F$ to be a pseudo-finite field not containing the square roots of $2$. Then $f$ further induces a bijection in some sufficiently big finite field not containing the square roots of 2. This however is not possible for cardinality reasons.

We typically need this kind of comparison allowed by our preorder in equation \eqref{eq.ge 1} of the proof of Theorem \ref{t.bound on sum of variations}.
\end{ex}
Note that in case the residue field has definable Skolem functions, the inclusion of Lemma \ref{l.<C+} is an equality, as shown by the following lemma. 

\begin{lem}\label{l.<skol} 
If the residue field $k$ of $K$ has definable Skolem functions, then
for any $m,n,p,r\in \NN$ and any $\emptyset$-definable set 
$X\subseteq K^n\times k^m\times \ZZ^r$,  
we have 
$\mathrm{I}\mathcal{C}_{\succcurlyeq 0}(X) = \mathrm{I}\mathcal{C}_+(X)$.
\end{lem}

\begin{proof}
The inclusion $ \mathrm{I}\mathcal{C}_+(X) \subseteq \mathrm{I}\mathcal{C}_{\succcurlyeq 0}(X)$  is Lemma~\ref{l.<C+}. 
Let now $F \in \mathrm{I}\mathcal{C}_{\succcurlyeq 0}(X)$ be given, witnessed by
$Z \overset{f}{\twoheadrightarrow} Y$, with $Y,Z\subset k^{m'}$ for $m'\in \NN$, and $\varphi \in \mathrm{I}\mathcal{C}_+(Y \times X)$. That is to say $F = \pi_{X!}(f_{X}^*\varphi) - \pi_{X!}(\varphi)$.

By the assumption that the residue field has definable Skolem functions (and using that the value group $\ZZ$ also has definable Skolem functions) we find a $\emptyset$-definable map $g\colon Y \to Z$ such that $f\circ g$ is the identity on $Y$.
Set $Z' := Z \setminus g(Y)$, and let $f'\colon Z' \to Y$ be the restriction of $f$ to $Z'$.
Then
$F =  \pi_{X!}(f'^{*}_{X}\varphi)$, and $\pi_{X!}(f'^{*}_{X}\varphi)$ lies in $\mathrm{I}\mathcal{C}_+(X)$, since $\varphi$ lies in $\mathrm{I}\mathcal{C}_+(Y \times X)$ (using in particular Remark \ref{r.Fubini}).
\end{proof}

\begin{lem}\label{l.<+} For $m,n,r\in \NN$ and $X\subseteq K^n\times k^m\times \ZZ^r$ a $\emptyset$-definable set, let $F_1, F_2\in \mathrm{I}\mathcal{C}_{\succcurlyeq 0}(X)$. Then 
$F_1+F_2\in \mathrm{I}\mathcal{C}_{\succcurlyeq 0}(X)$.

In particular, we obtain a partial preorder on $\mathcal{C}(X)$ by setting 
$$F \preccurlyeq G \Longleftrightarrow G - F \in \mathrm{I}\mathcal{C}_{\succcurlyeq0}(X).$$
\end{lem}

\private{
Question: What kind of equiv rel is defined by $F \le G \wedge G \le F$?
}
\begin{proof}
Suppose that   $F_1, F_2 \in \mathrm{I}\mathcal{C}_{\succcurlyeq0}(X)$ are witnessed by $Z_i \overset{f_i}{\longtwoheadrightarrow} Y_i$ and $\varphi_i \in \mathrm{I}\mathcal{C}_+(Y_i \times X)$, that is to say $F_i =  \pi_{X!}(f_{iX}^*\varphi_i) - \pi_{X!}(\varphi_i)$.
Without loss, we can assume that $Y_1$ and $Y_2$ are disjoint subsets in 
$k^m$, for some $m\in \NN$ and similarly for $Z_1$ and $Z_2$.
Set
$Y = Y_1 \sqcup Y_2$,
$Z = Z_1 \sqcup Z_2$, define $f\colon Z \to Y$ such that $f_{\restriction Z_i} = f_i$ and $\varphi \in \mathrm{I}\mathcal{C}_+(Y \times X)$ such that $\varphi_{\restriction Y_i \times X} = \varphi_i$. Then 
$F_1 + F_2 = \pi_{X!}(f_X^*\varphi) - \pi_{X!}(\varphi)$.
\end{proof}

\begin{rem}\label{rem.multiplication by a measure}
For $m,n,r\in \NN$ and $X\subseteq K^n\times k^m\times \ZZ^r$ a $\emptyset$-definable set, when $F\in \mathrm{I}\mathcal{C}_{\succcurlyeq 0}(X)$ and $G\in \mathcal C_+(X)$, then $F \cdot G\in \mathrm{I}\mathcal{C}_{\succcurlyeq 0}(X)$, since in 
Definition~\ref{d.<}, we can replace $\varphi$ by $\varphi\cdot G$.

\end{rem}

\begin{qu}
However, we do not know whether 
$F, G\in \mathrm{I}\mathcal{C}_{\succcurlyeq 0}(X)$ implies $F\cdot G\in \mathrm{I}\mathcal{C}_{\succcurlyeq 0}(X)$.
\end{qu}

The partial preorder on $\mathrm{I}\mathcal{C}(X)$ is compatible with pull-back.

\begin{lem}\label{l.<*} For $m,n,r\in \NN$ and  $\emptyset$-definable sets $X, \tilde{X}\subseteq K^n\times k^m\times \ZZ^r$, if 
  $g\colon \tilde X \to X$, $1_{\tilde{X}}$ is in $\mathrm{I}_X\mathcal{C}_+(\tilde{X}\overset{g}{\to} X)$ and $F \in \mathrm{I}\mathcal{C}_{\succcurlyeq0}(X)$, then  $g^*F \in \mathrm{I}\mathcal{C}_{\succcurlyeq0}(\tilde X)$.
\end{lem}

\private{
If $g$ is surjective, we should maybe even claim ``if and only if''. But how does one prove this?}

\begin{proof} 
Let $F \in \mathrm{I}\mathcal{C}_{\succcurlyeq 0}(X)$ be witnessed by
$Z \overset{f}{\longtwoheadrightarrow} Y$, and $\varphi \in \mathrm{I}\mathcal{C}_+(Y \times X)$. That is to say
$F = \pi_{X!}(f_{X}^*\varphi) - \pi_{X!}(\varphi)$.
Set $\tilde \varphi := g_Y^*\varphi$, where $g_Y\colon Y \times \tilde X \to Y \times X, (y,\tilde x) \mapsto (y, g(\tilde x))$.
Then $\tilde\varphi \in \mathrm{I}\mathcal C(Y \times \tilde X)$, and
$g^* F \in \mathrm{I}\mathcal{C}_{\succcurlyeq 0}(X)$ is witnessed by $Z \overset{f}{\longtwoheadrightarrow} Y$ and $\tilde \varphi$ (see
\cite[Theorem 10.1.1 A3]{CL.mot} for $\tilde \varphi\in \mathrm{I}\mathcal{C}_+(Y\times \tilde{X})$).
\end{proof}

The partial order is compatible with integration, in the following sense.

\begin{lem}\label{l.<int}
Let  $m,n,r\in \NN$, $U,V\subseteq K^n\times k^m\times \ZZ^r$ be $\emptyset$-definable sets, $F_1,F_2\in  \mathrm{I}_V\mathcal{C}(U\times V)$ and $\pi_V: U\times V\to V$ the standard projection. 
If $F_1 \preccurlyeq F_2$, then $\pi_{V!}F_1 \preccurlyeq \pi_{V!}F_2$.
This may be written, using the integral notation 
$$
\int_U F_1(u,\cdot) \ du \preccurlyeq \int_U F_2(u,\cdot) \ du. 
$$
\end{lem}

\begin{proof}
By additivity of the integral, it suffices to prove that for $F \in \mathrm{I}\mathcal{C}_{\succcurlyeq0}(U \times V)
$, we have
$\tilde F := \pi_{V!}F \in \mathrm{I}\mathcal{C}_{\succcurlyeq0}(V)$.

Let $Z \overset{f}{\longtwoheadrightarrow} Y$,   $\varphi \in  \mathrm{I} \mathcal{C}_+(Y \times U\times V)$
witness that $F \succcurlyeq 0$. That is to say $F = \pi_{U \times V!}(f_{U \times V}^*\varphi) - \pi_{U \times V !}(\varphi)$.

By Remark \ref{r.Fubini} 
$\varphi\in \mathrm{I}_{Y\times V}\mathcal{C}_+(Y\times U\times V)$. We can therefore define
$
\tilde \varphi := \pi_{Y \times V!}\varphi$, and again by Remark \ref{r.Fubini}, we have
$\tilde \varphi \in \mathrm{I}\mathcal{C}_+(Y \times V)$.

Then  $Z \overset{f}{\longtwoheadrightarrow} Y$ and $\tilde \varphi$ witness $\tilde F \ge 0$. 
Indeed, we have
\begin{align*}
\tilde F &=  \pi_{V!}\pi_{U\times V!}f_{U \times V}^*\varphi - \pi_{V!}\pi_{U\times V!}\varphi
\\
&= \pi_{V!}\pi_{Z \times V!}f_{U \times V}^*\varphi - \pi_{V!}\pi_{Y \times V!} \varphi
\\
&= \pi_{V!}f_{V}^*\pi_{Y \times V!}\varphi - \pi_{V!}\pi_{Y \times V!} \varphi
\\
&=\pi_{V!}f_{V}^*\tilde \varphi - \pi_{V!}\tilde \varphi.
\end{align*}
Note that in the comptation above $\pi_V$ does not always represent the same projection on $V$, but nevertheless no confusion is possible.

\end{proof}

For our partial preorder to work nicely in families, one might want, in Definition~\ref{d.<}, to replace $Y \times X$ by some set $V \subset X \times k^{m'}$ where the fibre $V_x$ may depend on $x$, and similarly for $Z$ and $f$. The following lemma says that this is indeed possible.

\begin{lem}\label{l.<surj}
Let $X\subseteq K^n\times k^m\times \ZZ^r$,
$V \subset X \times k^{m'}$, $W \subset X \times k^{m''}$
be $\emptyset$-definable sets, for some
$m,n,r,m',m''\in \NN$, 
let $g\colon W \to V$ be a surjective $\emptyset$-definable map, and let $\psi \in \mathrm{I}\mathcal{C}_+(V)$ be a constructible motivic function. We denote by 
$\pi_X^V:V\to X$ and $\pi_X^W:W\to X$ the standard projections, and for $x\in X$, by $V_x$ and $W_x$ their fibre over $x$.
We assume that $g$ is fibrewise (over $X$), that is, the following diagram is commutative, 
$$
  \xymatrixcolsep{2pc}\xymatrix{
   W \ar@{->>}[dr]_{\pi^W_X} \ar@{->>}[rr]^{g}  &   & 
   V \ar@{->>}[ld]^{\pi^V_X}   &  \\
 & X  & 
  } 
$$
and we assume moreover that $g$ has finite fibres. Then  
 \[F := \pi_{X!}^{W}(g^*\psi) - \pi_{X!}^{V}(\psi) \in \mathrm{I}\mathcal{C}_{\succcurlyeq0}(X).
 \]
\end{lem}

\begin{proof}
By Remark~\ref{r.k.st.e}, applied to $V$, $W$ and to the graph of $g$, there exist $\hat m \in \NN$, a definable map $\rho\colon X \to k^{\hat m}$, and $\emptyset$-definable $Y \subset k^{\hat m} \times k^{m'}$, $Z \subset k^{\hat m} \times k^{m''}$, $f\colon Z \twoheadrightarrow Y$ such that for every $x \in X$, we have $V_x = Y_{\rho(x)}$, $W_x = Z_{\rho(x)}$ and $g|_{W_x} = f|_{Z_{\rho(x)}}$.
Define $h\colon V \to X\times Y$ by $h(x,v) = (x, \rho(x), v)$ (where $x \in X$ and $v \in k^{m'}$) and
$\varphi := h_!\psi \in \mathrm{I}\mathcal{C}_{\succcurlyeq0}(X \times Y)$. We claim that $f\colon Z \twoheadrightarrow Y$ and $\varphi$ witness that $F$ lies in $\mathrm{I}\mathcal{C}_{\succcurlyeq0}(X)$, i.e. that
$F = \pi_{X!}f_X^*\varphi - \pi_{X!}\varphi$.
More precisely, we claim that 
\begin{itemize}
\item[(a)] 
$\pi_{X!}f_X^*\varphi = \pi_{X!}^{W}(g^*\psi)$, and 
\item[(b)]
$\pi_{X!}\varphi = \pi_{X!}^{V}(\psi)$.
\end{itemize}

(b) follows directly from $\pi_X \circ h = \pi^V_X$ (namely, $\pi_{X!}\varphi = \pi_{X!}h_!\psi= \pi_{X!}^{V}(\psi)$).

For (a), define $\tilde h\colon W \to X\times Z$ by $\tilde h(x,w) = (x, \rho(x), w)$ (where $x \in X$ and $w \in k^{m''}$), and note that we have $f\circ \tilde h = h \circ g$. Then we have
\[
\pi_{X!}f_X^*\varphi
=\pi_{X!}f_X^*h_!\psi
=\pi_{X!}\tilde h_!g^*\psi
=\pi_{X!}^{W}(g^*\psi),
\]
where the middle equality becomes clear if we use $h$ and $\tilde h$ to identify $V$ with $h(V) \subset X \times Y$ and $W$ with $\tilde h(W) \subset X \times W$: With this identification, we have that
$g$ is the restriction of $f_X$ to $h(W) \to h(V)$,
$h_!\psi$ is the extension by $0$ of $\psi\in \mathrm{I}\mathcal{C}_{\succcurlyeq0}(h(V))$ to $X \times Y$, and $\tilde h_!g^*\psi$ is the extension by $0$ of $g^*\psi\in \mathrm{I}\mathcal{C}_{\succcurlyeq0}(\tilde h(W))$ to $X \times Z$.
\end{proof}

\begin{rem}\label{r.<sorts}
In Lemma~\ref{l.<surj}, one could more generally allow
$V \subset X \times k^{m'} \times \ZZ^{r'}$ (and similarly for $W$), and even $V \subset X \times K^{n'} \times k^{m'} \times \ZZ^{r'}$ if we use the $0$-dimensional motivic measure on $K^{n'}$. Indeed, one can still find $\rho\colon X \to k^{\hat m}$ as in the proof of the lemma, using orthogonality of $k$ and $\ZZ$, and using that finite subsets of $K^{n'}$ are in definable bijection with definable subsets of $k^{m'''}$.
\end{rem}

\section{Motivic Vitushkin invariants}\label{s.Motivic Vitushkin invariants}

In the following, $X$ is always a bounded $\emptyset$-definable subset of $K^n$. For the moment, we assume $n \ge 1$.

\begin{defn}\label{d.V0}
Let $B_1, \dots, B_\ell \subseteq K^n$ be the balls and singletons provided by Theorem~\ref{t.S0}, i.e. such that $X$ is not riso-trivial on a ball $B$ if and only if $B$ contains one of the $B_i$, and let $S_{0}$ be a finite $\emptyset$-definable subset of $k^m$, for some integer $m$, definably parameterizing the $B_i$, $i=1, \ldots, \ell$, i.e. such that there exists a $\emptyset$-definable surjection
$\bigcup_i B_i \to S_{0}$ whose fibres are exactly the $B_i$. Note that $S_{0}$ is uniquely determined up to $\emptyset$-definable bijection. 

We then define the \emph{$0$-dimensional motivic Vitushkin variation}, corresponding to the number of connected components in the real setting, as
$$
V_0(X) := \mu_0(S_{0}) \in \mathcal{C}_+(\{pt\}),$$ 
where $\mu_0$ denotes the $0$-dimensional motivic measure (or motivic cardinality) of Section \ref{s.Motivic integrals}. 

For subsets of $K^0$, we set $V_0(\emptyset) := 0$ and $V_0(K^0) := 1$.

\end{defn}

\begin{rem}\label{r.V0fin}
For finite sets $X \subseteq K^n$, we simply have $V_0(X) = \mu_0(X)$. Indeed, in this case,
$B_1, \dots, B_\ell$ are just the singletons of $X$ and thus $S_{0}$ is in $\emptyset$-definable bijection with $X$, provided that $n \ge 1$. In the case $n = 0$, we intentionally defined
$V_0$ to make this remark hold.
\end{rem}

If, instead of a single $\emptyset$-definable set $X$, we have a $\emptyset$-definable family $(X_z)_{z \in Z}$, then one can naturally adapt
Definition~\ref{d.V0} to yield
a constructible motivic function $(V_0(X_z))_{z \in Z} \in \mathcal \mathcal{C}_+(Z)$. The same also applies to the definition of $V_i(X)$ below, and also to the relative versions $V_i(X, B)$, where both, $X$ and $B$ can vary in families.

\begin{rem}\label{r.V0.bound}
For any $\emptyset$-definable family $(X_z)_{z \in Z}$ of sets $X_z \subset K^n$, the constructible motivic function $(V_0(X_z))_z$ is uniformly bounded by a constant function in the sense of Definition~\ref{d.convergence Mot}.
Indeed, by Remark~\ref{r.S0.bound}, the number $\ell$ of balls $B_i$ is uniformly bounded, which implies that they can be parametrized by a subset of $k^m$ for $m$ not depending on $z$. Thus,
applying Remark \ref{r.subset and inequality},
the constructible motivic function on $Z$ which is constant equal to $[k^m] = \LL^m$ is a possible bound.
\end{rem}

\medskip

\begin{notn}\label{n.normalisation measure}
We denote by $\mathrm{d} g$ the  $\GL_n(\OK)$-invariant motivic measure on $\GL_n(\OK)$ (see \cite[Section 15]{CL.mot} for the notion of motivic measure on manifold), normalized such that the total measure of $\GL_n(\OK)$ is equal to $[\GL_n(k)]\cdot \LL^{-\dim (\GL_n)}$.
\end{notn}

Fix $i\in \{0, \ldots, n\}$. There is a natural $K^n \rtimes \GL_n(\OK)$-invariant motivic measure on the Grassmanian manifold $\affgr{n-i}$ of $(n-i)$-dimensional affine subspaces of $K^n$. Formally, it can be defined using affine charts, but computations become more natural with the following definition, which is equivalent (up to scaling).

\begin{notn}\label{n.mesure grassmannienne}
Given an affine space $\bar{P} \in \affgr{n-i}$, we write $P \in  \vecgr{n-i}$ for the direction of $\bar P$.

We then define a measure on $ \affgr{n-i}$. To this end, we fix once and for all a $\emptyset$-definable element $g_0\in \GL_n(\OK)$ (e.g. the identity matrix) and set
$P_0:=g_0\cdot(K^{n-i}\times \{0\}^{i})$.

Those choices yield a surjective map $\gamma\colon \GL_n(\OK) \to  \vecgr{n-i}$ sending $g$ to $g\cdot P_0$, and a surjective map $\bar\gamma \colon \GL_n(\OK)\times K^i \to  \affgr{n-i}$,
sending $(g, y)$ to $gg_0\cdot( K^{n-i}\times \{y\})$. Clearly, those two maps are compatible with respect to taking the direction of an affine subspace of dimension $n-i$ of $K^n$, and note that once we fix $g \in \GL_n(\OK)$, $\bar\gamma$ allows us to identify $K^i$ with $K^n/(g\cdot P_0)$. We denote by $\pi_g\colon K^n \to K^i$ the projection from $K^n$ to $K^n/(g\cdot P_0)$ using that identification. More accurately, 
we have for $x\in K^n$ (or for any element of $K^n$ in the same class as $x$ modulo $g\cdot P_0$), $\pi_g(x)=y$, where $(gg_0)^{-1}(x)=a+y$ with $a\in K^{n-i}\times \{0\}^i$ and $y\in \{0\}^{n-i}\times K^i$.
In particular, for $y\in K^i$ and $g\in \GL_n(\OK)$, we have $\pi_g^{-1}(y)=\bar\gamma(g,y)$. 

We use $\gamma$ to push forward the measure on $\GL_n(\OK)$ to $ \vecgr{n-i}$, i.e. for any $\varphi\in \mathcal{C}_+( \vecgr{n-i})$, we say that $\varphi$ is integrable if $\varphi\circ \gamma\in \mathrm{I}\mathcal{C}_+(\GL_n(\OK))$ with respect to the motivic measure $\mathrm{d}g$, and we define the integral of $\varphi$ by
\begin{align*}
  \int_{P\in  \vecgr{n-i} } \varphi(P) \ \mathrm{d}P & :=
  \int_{g\in \GL_n(\OK)} \varphi\circ \gamma(g) \ \mathrm{d}g
\end{align*}
In a similar way, we use $\bar\gamma$ to push forward the motivic product measure
$dg \otimes dy$ on $\GL_n(\OK)\times K^i$ to $ \affgr{n-i}$. In particular, for $\varphi \in \mathrm{I}\mathcal{C}_+(
 \affgr{n-i})$, we set
\begin{align*}
 \int_{\bar{P}\in  \affgr{n-i}} \varphi(\bar{P}) \ \mathrm{d}\bar{P} & :=
  \int_{g\in \GL_n(\OK)} \int_{y\in K^i} \varphi\circ  \bar{\gamma}(g,y) \ \mathrm{d}y  \ \mathrm{d}g \\
&  =  \int_{y\in K^i} \int_{g\in \GL_n(\OK)}  \varphi\circ \bar{\gamma}(g,y) \ \mathrm{d}g  \ \mathrm{d}y. 
 \end{align*}
 \end{notn}
 
\begin{rem}\label{r.identifying Ki and Kn/P}
We note, with this definition of integration on 
$ \affgr{n-i}$, that the  choice of $g_0\in \GL_n(\OK)$ and $P_0$ made to identify $K^n/P$ and $K^i$
does not matter by the change of variable theorem.
So let us, for the remainder of this paper, set $g_0=I_n$ (and thus $P_0=K^{n-i}\times \{0\}^i$). 
\end{rem}

\begin{defn}[Multidimensional motivic Vitushkin's variations]\label{d.Vi} For $X$ a $\emptyset$-definable subset of $K^n$ and for $i=1, \ldots, n$, we define $V_i(X) \in \mathcal C_+(\{pt\})$, the \emph{$i$th motivic Vitushkin variation of $X$}, by 
$$
V_i(X) := \int_{\bar P \in \affgr{n-i}} V_0(X \cap \bar P) \ \mathrm{d}\bar P, 
$$
\end{defn}

\begin{rem}
For a $\emptyset$-definable bounded subset $X$ of $K^n$ and a general $\bar{P}\in \affgr{n-i}(k)$, the set $X\cap \bar{P}$ is not $\emptyset$-definable. 
Nevertheless the subset $Y:=\{(x,\bar{P}); x\in X\cap \bar{P}, \bar{P}\in \affgr{n-i}\}$ of $X\times \affgr{n-i}$ is $\emptyset$-definable, since $\affgr{n-i}$ is itself $\emptyset$-definable. 
Denoting $\pi: Y \to \affgr{n-i}$ the projection map $\pi(x,\bar{P})=\bar{P}$, in view of Remark \ref{r.notation famille},
the $\emptyset$-definable family $ (V_0(X \cap \bar P))_{\bar{P}\in \affgr{n-i}(k)}$
 stands for $\pi_{!\affgr{n-i}}(\mathbf{1}_Y)\in 
 \mathcal{C}_+( \affgr{n-i})$, according to Remark \ref{rem.direct product}. Consequently, $V_i(X)$, the motivic integral of this constructible motivic function on $\affgr{n-i}$, makes sense.

Note also that since $X$ is bounded, 
$\bar P \mapsto V_0(X \cap \bar P)$ is non-zero only on a bounded subset of $\affgr{n-i}$, so that the integral $V_i(X)$ is finite.
\end{rem}

\begin{rem}
\label{r.Vi}
The following properties of $V_i$ follow from the nature of the measure 
considered on $\affgr{n-i}$. 
\begin{enumerate}
 \item
 \label{item.homogeneity} 
 For any $i=1, \ldots, n$, for any 
 $r \in K^\times$, $$V_i(r\cdot X) =|r|^i V_i(X).$$
\item For any  $i=1, \ldots, n$ and any $j\le i$, 
$$V_i(X)=\int_{\bar{P}\in  \affgr{n-j} } V_{i-j}(X\cap \bar{P})\ \text{d}\bar{P}.$$
\end{enumerate}
\end{rem}
\private{Is the following true? Invariance of $V_0, V_i$ under definable risometry... même non-definable risometry?
}

\begin{rem}
Note that in contrast to the real setting, our Vitushkin invariants are not sub-additive. In the real case the sub-additivity of the $V_i$'s comes from the sub-additivity of $V_0$ itself. But in the nonarchimedean case, $V_0$ is not sub-additive. Intuitively, the reason is that it rather counts singularities than connected components.
As an example, consider $X_1 := \OK \times \{0\}$, and let $X_2$ be the graph of the function $\OK \to \OK, x \mapsto t(x^3 - x)$, where $t \in \OK$ has valuation $1$. Then for both, $i = 1,2$, Theorem~\ref{t.S0} yields only the ball $t^{-1}\OK \times t^{-1}\OK$, so we have $V_0(X_1) = V_0(X_2) = 1$.
On the other hand, $(-1,0), (0,0),(1,0)$ are singularities of $X_1 \cup X_2$, so that we obtain
$V_0(X_1 \cup X_2) = \mu_0(\{(-1,0),(0,0),(1,0)\})=3$.
\end{rem}

Note however that by the Cauchy-Crofton formula (see Theorem~\ref{t.Crofton}), we do have sub-additivity (even additivity of $V_i(X)$ for sets $X$ of dimension~$i$.)

\section{Motivic Cauchy-Crofton formula}\label{s.Motivic Cauchy-Crofton}

We now show that in case $d=\dim(X)$, $V_d(X)$ is nothing else than the motivic measure of $X$ (up to multiplication by some universal constant); this is the motivic version of the classical Cauchy-Crofton formula, (a local version of the Cauchy-Crofton formula involving the local density has been proved in the real case in  \cite{Co1}, in the $p$-adic case in \cite{CCL.cones}, and in  \cite{For.LocalCauchy} in the general nonarchimedean setting). 
This may be an indication of the relevance of 
our nonarchimedean definition of $V_0$ in terms of riso-triviality as a substitute of the real counting points measure involved in the real Cauchy-Crofton formula.  

\begin{thm}[Motivic Cauchy-Crofton formula]\label{t.Crofton} Let $X$ a bounded $\emptyset$-definable subset of 
$K^n$ of dimension $d$. Then $V_i(X)=0$ for $i>d$ and
\begin{equation}\label{eq.CC}
V_d(X) = C(n,d)\mu(X),
\end{equation}
in $\mathcal{C}_+(\{pt\})$,
for some universal constant $C(n,d)\in \mathcal{C}_+(\{pt\})$ depending only on $d$ and $n$, and explicitly given in Lemma \ref{l.Crofton lineaire}.

Moreover, the same formula 
also holds in families, that is, given a $\emptyset$-definable set $Z\subseteq K^m$
and a $\emptyset$-definable family
$(X_z)_{z \in Z}$ of sets $X_z\subseteq K^n$ of dimension $d$, with $X_z$ uniformly bounded, 
we have 
$$\left( V_d(X_z) \right)_z= (C(n,d) \mu(X_z))_z $$
in $\mathcal C_+(Z)$.
\end{thm}

We first give the proof of Theorem \ref{t.Crofton} in the simple case where 
$X$ is a $\emptyset$-definable subset of dimension $d$ of a $d$-dimensional affine space of $K^n$. 

\begin{lem}\label{l.Crofton lineaire}
The motivic Cauchy-Crofton formula 
$$
V_d(X) = C(n,d)\mu(X),
$$
of Theorem \ref{t.Crofton} holds for 
$X$ a $\emptyset$-definable subset of dimension $d$ of a $d$-dimensional affine space $V$  of $K^n$, where
$$C(n,d) = \int_{g\in \GL_n(\OK)}  \LL^{\mathrm{ordjac}\pi_{g\restriction V}}
 \ \mathrm{d}g,$$
 and for $g\in \GL(\OK)$, $\pi_g: K^n\to K^d$ is the projection from $K^n$ to $K^n/(g\cdot P_0)$, $K^n/(g\cdot P_0)$ being identified to $K^d$ (see Notation \ref{n.mesure grassmannienne} and
 Remark \ref{r.identifying Ki and Kn/P}).  
 In particular one can take $C(n,d)=1$ using the convenient normalisation of the measure $\mathrm{d}g$ on $\GL_n(\OK)$. 
 
 As mentioned in Theorem \ref{t.Crofton}, this formula also holds in families. 
\end{lem}

\begin{proof}[Proof of Lemma \ref{l.Crofton lineaire}]
We have by definition of the measure $\mathrm{d}\bar{P}$ 
\begin{align*}
V_d(X)= & \int_{\bar{P}\in \affgr{n-d}} V_0(X\cap \bar{P}) \ \mathrm{d}\bar{P} \\
&= \int_{g\in \GL_n(\OK) } \int_{y\in K^d} V_0(X\cap \pi_g^{-1}(y)) \ \mathrm{d}y \ \mathrm{d}g. 
\end{align*}
For almost all $g\in \GL_n(\OK)$, the restriction 
$\pi_{g\restriction V}$ of $\pi_g$ to $V$ is a bijection from $V$ to $K^d$. By the change of variable formula, we thus have 
$$
\int_{y\in K^d} V_0(X\cap \pi_g^{-1}(y)) \ \mathrm{d}y=
\int_{v\in V} \LL^{\mathrm{ordjac}\pi_{g\restriction V}(v)} \mathbf{1}_X(v) \ \mathrm{d}v.
$$
Noting that $\mathrm{ordjac}\pi_{g\restriction V}(v)$ does not depend on $v\in V$, and denoting $\mathrm{ordjac}\pi_{g\restriction V}$ this constant, we get 
\begin{align*}
V_d(X)= & \int_{g\in \GL_n(\OK)}  \LL^{\mathrm{ordjac}\pi_{g\restriction V}}
\int_{v\in V}  \mathbf{1}_X(v) \ \mathrm{d}v
\ \mathrm{d}g \\
 =& \mu(X) \int_{g\in \GL_n(\OK)}  \LL^{\mathrm{ordjac}\pi_{g\restriction V}}
 \ \mathrm{d}g.
\end{align*}
We denote by $C(n,d)$ the constant  $\int_{g\in \GL_n(\OK)}  \LL^{\mathrm{ordjac}\pi_{g\restriction V}}
 \ \mathrm{d}g$ and observe that $C(n,d)$ does not depend on the choice of $V$, since $\mathrm{d}g$ is 
$\GL_n(\OK)$-invariant.
\end{proof}

Before proving Theorem~\ref{t.Crofton} for general $X$, we prove an additivity result that will be needed in the proof.
\begin{lem}\label{l.addit}
Suppose that $X \subset K^n$ is an  $\emptyset$-definable set of dimension $d$ and that $X = \bigcup_{z \in Z} X_z$ is a partition of $X$ into subsets $X_z$ such that $(X_z)_{z \in Z}$ is a $\emptyset$-definable family, where $Z \subset k^m \times \ZZ^r$ for some $m$ and $r$.
Then $$
V_d(X) = \int_Z V_d(X_z)
$$
\end{lem}

\begin{proof}
For generic $\bar P$, the set $X \cap \bar P$ is a finite number of points, so that $V_0(X \cap \bar P) = \mu(X \cap \bar P)$ and $V_0(X_z\cap \bar P) = \mu(X_z\cap \bar P)$ by Remark~\ref{r.V0fin}.
Thus, for generic $\bar P$, we have
$$\int_{z \in Z}  V_0(X_z\cap \bar{P}) \ \mathrm{d}x
= V_0(X \cap \bar{P}),$$  and hence, using Fubini theorem,
\begin{align*}
V_d(X)= & \int_{\bar{P}\in \affgr{n-i}} V_0(X\cap \bar{P}) \ \mathrm{d}\bar{P} \\
&= \int_{z \in Z} \int_{\bar{P}\in \affgr{n-i}} V_0(X_z\cap \bar{P}) \ \mathrm{d}\bar{P} 
\\
&=\int_{z \in Z} V_d(X_z).\qedhere
\end{align*}
\end{proof}

In the proof of Theorem~\ref{t.Crofton}, it will be handy to have a metric (taking values in $\LL^\ZZ \cup \{0\}$) on the Grassmannians 
$\vecgr{d}$. A natural metric is the following one:
\begin{defn}\label{d.Delta}
For $P_1,P_2 \in \vecgr{d}$, we set
\[
\Delta(P_1, P_2) := \max\{\min \{|a_1-a_2|; a_1 \in P_1\}; a_2 \in P_2, |a_2| \le 1\}.
\]
\end{defn}
\begin{rem}\label{r.tub.grass}
Note that Remark~\ref{r.tub} about the measure of tubular neighbourhoods also applies in $\vecgr{d}$, using this metric and the measure from Notation~\ref{n.mesure grassmannienne}, namely: Given a $d$-dimensional 
$\emptyset$-definable set $X \subset \vecgr{d}$ and some $\nu \in \LL^\ZZ$, the tubular neighbourhood  
\[
\{P \in \vecgr{d};\exists P'\in X, \Delta(P-P')\le \nu\}
\]
has measure bounded by $C\cdot \nu^{N-d}$ for some constant $C \in \mathcal{C}_+(\{pt\})$ not depending on $\nu$, and where $N := \dim( \vecgr{n-d})= d(n-d)$.
Indeed, we can cover $\vecgr{d}$ by finitely many subsets of $K^N$ in such a way that the measures and the metrics from $\vecgr{d}$ and $K^N$ differ by at most some factor.
\end{rem}

We now consider the general case where $X$ is a $\emptyset$-definable bounded subset of $K^n$ of dimension $d$ to prove Theorem \ref{t.Crofton}. The proof of the family version is a straight forward adaptation of the non-family version, so for the sake of readability, we only treat the absolute case.
\begin{proof}[Proof of Theorem \ref{t.Crofton}]
If $i>d$, we have $\dim(\pi_g(X))\le d<i$ and therefore for any 
$y$ not in the set $\pi_g(X)$ of measure zero, $\bar{P}=\pi_g^{-1}(y)$ does not encounter $X$. For such generic $\bar{P}$, $V_0(X\cap \bar{P})=0$, which proves that $V_i(X)=0$. This also shows that in the case $i = d$, which we prove now, we can freely remove ($\emptyset$-definable) subsets of positive codimension from $X$.

\medskip
{\bf Step 1.}\refstepcounter{dummy}\label{step1}
We first apply \cite[Theorem 5.7.3]{iCR.hmin} together with Addendum~5 of \cite[Theorem 5.2.4]{iCR.hmin} to get a decomposition of $X$ into sets which are, up to permutation of coordinates, reparametrized cells with 1-Lipschitz centres.
More precisely
(using also Remarks~5.3.2 and 5.7.2 of \cite{iCR.hmin}), we obtain a partition of $X$ into  sets $X_\xi$, where $\xi$ runs over $k^m$ for some integer $m$, such that 
$(X_\xi)_{\xi \in k^m}$
is a $\emptyset$-definable family
and each $X_\xi$ is, up to a permutation of coordinates, the graph of a $1$-Lipschitz map $f_\xi\colon \bar X_\xi \subseteq K^d \to K^{n-d}$ (where $\bar X_\xi$ is the projection of $X_\xi$ to $K^d$).
Now note that the theorem for $X$ follows from the family version of the theorem for $(X_\xi)_\xi$,
since
$$
\mu(X) = \int \mu(X_\xi)d\xi
$$
and
$$
V_d(X) = \int V_d(X_\xi)d\xi
$$
by Lemma~\ref{l.addit}.
We will therefore assume that $X$ itself is the graph of a $1$-Lipschitz map $f\colon \bar X \subseteq K^d \to K^{n-d}$ (and let the reader adapt that proof to a family version).
We may moreover assume in this process that the dimension of $X$ did not drop, i.e. $\dim X = \dim \bar X = d= i$.

We now control how well $f$ is approximated by its degree 1 Taylor polynomial on balls $B \subset \bar X$. To this end, we will remove various $\emptyset$-definable subsets of $\bar X$ of positive co-dimension, namely as follows.
Applying \cite[Theorem 5.6.1]{iCR.hmin} to each coordinate function of $f$ (after extending $f$ trivially to all $K^d$) yields that, for every ball $B \subseteq \bar X$, we have the following
\begin{itemize}
    \item[-] $f_{\restriction B}$ is $C^2$,
    \item[-] The norms
      $|\partial^2 f/\partial x_{i} \partial x_j|$, $1 \le i, j \le d$, are constant on $B$,
    \item[-] For every $x, x' \in B$, we have
    \begin{equation}\label{eq.5.6.1}|f(x) - T_{f,x'}^{\le 1}(x)| \le \max_{i,j \le d}|\partial^2 f/\partial x_{i} \partial x_j|\cdot |x_i-x'_i|\cdot |x_j-x'_j|,\end{equation}
    where $T^{\le 1}_{f,x'}(x)=f(x')+Df_{(x')}\cdot (x-x')$ is the first order Taylor polynomial in $x$ of $f$ at $x'$.
\end{itemize}

Next, for $a\in B$ and $j\in \{1, \ldots, n\}$, let  
$$
g_a:x_j\mapsto  \partial f/\partial x_i(a_1, \ldots, a_{j-1}, x_j, a_{j+1}, \ldots, a_n).
$$
For each $j$ and $a$, apply the proof of \cite[Corollary~3.2.5]{iCR.hmin} with $r=0$ to $g_a$ and remove the bad points from $\bar X$, so that for $B \subseteq \bar X$ and $x_0 \in B$ we additionally have the inequality
\begin{equation}\label{eq.taylor}
|g_a'(x_0)|\cdot \rad B < |g_a(x_0)|
\end{equation}
obtained in that proof. Note that the open radius of $B$ used in \cite{iCR.hmin} is equal to $\LL\cdot \rad B$ (using our notation), which is why we obtain a strict inequality in
\eqref{eq.taylor}.
Using that $f$ is assumed to be $1$-Lipschitz, \eqref{eq.taylor} implies
$$
|\partial^2 f/\partial x_{i} \partial x_j| =
|g_a'| < 1/\rad B.
$$
Plugging this bound into (\ref{eq.5.6.1}) yields
\begin{equation}\label{eq.Taylor}
 |f(x) - T_{f,x'}^{\le 1}(x)| < |x-x'|^2/\rad B.   
\end{equation}

\medskip
{\bf Step 2.} 
Given any $\lambda\in \LL^{\ZZ}$, let $\bar X_\lambda$ be the union of all balls $B$ of radius $\lambda$ contained in $\bar{X}$, and let $X_\lambda\subseteq X$ be the graph of the restriction of $f$ to $\bar X_\lambda$. Then $X \setminus \bigcup_{\lambda \in \LL^{\ZZ}}X_\lambda$ has positive codimension
and the sequence $(X_\lambda)_{\lambda\in \LL^{\ZZ}}$ is increasing for $\lambda \to 0$.

For any fixed $\lambda\in \LL^{\ZZ}$, we shall prove
the Cauchy-Crofton formula for $X_\lambda$, and then let $\lambda$ goes to $0$ to obtain the formula for $X$. Accordingly, from now on, we fix $\lambda\in \LL^\ZZ$ (formally, in all of the following arguments, we work in families, where $\lambda$ is a family parameter). 

Let $\nu\in \LL^\ZZ$ with $\nu \le \lambda$ and let
$B$ be a ball of radius $\nu$ 
such that $ X_\lambda \cap B \ne \emptyset$. (The idea is to show that the Cauchy-Crofton formula holds approximatively for $X_\lambda \cap B$, and that we can get arbitrarily good approximations by choosing $\nu$ small compared to $\lambda$).
Up to some harmless linear change of coordinates on $K^n$, we can assume that $0\in X_\lambda\cap B$ and that $T_0X=K^{d}\times \{0\}^{n-d}$.

Consider an affine plane $\bar{P}\in  \affgr{n-d}$ intersecting $X_\lambda\cap B$ at some point $(a,b)\in T_0X\oplus (\{0\} \times K^{n-d}) = K^n$ and intersecting $T_0X$ at some point~$c$. 

\medskip
{\bf Step 3.}
We now want to bound from above the measure of
the set of those $\bar{P} \in \affgr{n-d}$ (intersecting $X_\lambda \cap B$, as above), for which $c\notin B\cap T_0X$.
We start by bounding from above the measure of
the directions in $ \vecgr{n-d}$ of such $\bar{P}$.

Let us denote by $\mathcal{P}$ the variety of vector planes $P\in  \vecgr{n-d}$
not transverse to $T_0X$, that is the set of  $P\in  \vecgr{n-d}$ containing at least a line of $T_0X$. Since the flag variety of $P\in  \vecgr{n-d}$ containing a fixed line is of dimension $d(n-d)-d$, we have $\dim(\mathcal{P})=
d(n-d)-1=\dim( \vecgr{n-d})-1$. 
 
Since $(a,b)\in B$ and $T_{f,0}^{\le 1}=0$,
 by the Taylor approximation \eqref{eq.Taylor} one has $\vert b\vert 
 \le  \nu^2/\lambda$. On the other hand, since we consider points $c\in T_0X\setminus B$, we have
 $\vert a-c \vert >\nu$. 
 It follows that the direction of any line joining the point $(a,b)$ to a point $c\in T_0X\setminus B$ is in the $\nu/\lambda$-neighbourhood of a line of $T_0X$.
 As a consequence, any affine plane $\bar{P}\in  \affgr{n-d}$ containing 
$(a,b)$ and $c\in T_0X\setminus B$ has its direction $P$ in the $\nu/\lambda$-neighbourhood $\mathcal{P}_{\nu/\lambda}$ of $\mathcal{P}$.
By Remark~\ref{r.tub.grass}, the measure $\mu (\mathcal{P}_{\nu/\lambda})$ is bounded by $C\cdot(\nu/\lambda)^{\dim( \vecgr{n-d})-\dim(\mathcal{P})}=C\nu/\lambda$, for some constant $C \in \mathcal{C}_+(\{pt\})$.

Finally, let $\bar{\mathcal{P}}_{\nu/\lambda}$ be the set of $\bar{P}$ with direction in $\mathcal{P}_{\nu/\lambda}$ and $\bar{P} \cap B \ne \emptyset$.
Any $\bar{P}$ intersecting $X_\lambda\cap B$ but not intersecting $T_0X \cap B$ lies in this $\bar{\mathcal{P}}_{\nu/\lambda}$, so it remains to bound the measure of that set. To this end, fix a direction $P \in \mathcal{P}_{\nu/\lambda}$ and choose any 
$g \in \GL_n(\OK)$ satisfying $P = gP_0$. Then,
using Notation \ref{n.mesure grassmannienne}, the set of parameters $y\in K^n/P=K^d$ such that $\bar{P}=\pi_g^{-1}(y)$ and $\bar{P}$ intersects $B$, lies in $\pi_g(B)$, a set of measure $\nu^d$. It follows that the
preimage of $P$ under the direction map $\bar{\mathcal{P}}_{\nu/\lambda} \to \vecgr{n-d}$ has measure at most $\nu^d$, and hence $\bar{\mathcal{P}}_{\nu/\lambda}$ has measure at most
$C\nu^{d+1}/\lambda$.

More formally, by $\bar{\mathcal{P}}_{\nu/\lambda}$ having measure at most $C\nu^{d+1}/\lambda$, we mean that it is bounded by $C\nu^{d+1}/\lambda$ in the sense of Definition~\ref{d.convergence Mot}, considered as an element of $\mathcal{C}_+(Z)$, for some $\emptyset$-definable set $Z$ parametrizing our choices from Step~2.

\medskip
{\bf Step 4.}
Fix $P \in \vecgr{n-d} \setminus \mathcal{P}_{\nu/\lambda}$.
Then, for any $(a,b) \in X_\lambda\cap B$, there is a unique
$\bar P \in \bar{\mathcal{Q}}_{\nu/\lambda}$ intersecting
$X_\lambda\cap B$ in that point, and by Step 3, this 
$\bar{P}$ intersects $B\cap (K^d\times \{0\}^{n-d})$ at a point $c$.
We use this to define a map $\Pi$ (depending on $P$) from $B\cap (K^d\times \{0\}^{n-d})$ to $B\cap (K^d\times \{0\}^{n-d})$ sending $(a,0)$ to $c$. (Note that $b$ is determined by $a$ since $X_\lambda\cap B$ is the graph of a function.)

We now claim that $\Pi$ is a risometry (see Definition~\ref{d.riso}), that is to say  
$$\val \left( (\Pi(a) - \Pi(a')) - (a - a') \right) > \val(a - a'), \  \forall a,a' \in  B\cap (K^d\times \{0\}^{n-d}).$$
Indeed, on one hand, by \eqref{eq.Taylor}, $|b - b'|<  \nu/\lambda\cdot |a - a'|$.
And on the other hand, we have $c = a + Ab$, with $A \in \mathrm{Mat}_{d\times (n-d)}(K)$ with $|A| \le \frac{1}{\nu/\lambda}$.
Therefore 
$$|(\Pi(a) - \Pi(a')) - (a - a')| = |A(b-b')|
< \dfrac{1}{\nu/\lambda} \cdot \nu/\lambda\cdot |a - a'| = |a - a'|.$$

By Remark~\ref{r.riso.surj}, a risometry from a ball to itself is surjective.
This shows the converse of Step 3, that is, each 
$\bar{P} \in  \bar{\mathcal{Q}}_{\nu/\lambda}$ intersecting $T_0 X \cap B$ also intersects $X \cap B$.

We have defined in this way a definable family of bijections 
$$ (\varphi_P)_{\bar{P}\in\bar{\mathcal{Q}}_{\nu/\lambda}}$$ 
between 
$X_\lambda\cap B$ and $T_0X\cap B$, sending 
$X_\lambda\cap B\cap \bar{P}$ to $T_0X\cap B \cap \bar{P}$,
and therefore giving 
 \begin{equation}\label{eq.bijection}
\int_{\bar{P}\in \bar{\mathcal{Q}}_{\nu/\lambda}} V_0(T_0X\cap B\cap \bar{P}) \ 
\mathrm{d}\bar{P}
=
\int_{\bar{P}\in \bar{\mathcal{Q}}_{\nu/\lambda}} V_0(X_\lambda \cap B\cap \bar{P}) \ 
\mathrm{d}\bar{P}.
 \end{equation}  

\medskip
{\bf Step 5.}
We now stop assuming that our ball contains $0$ and that $T_0X = K^d \times \{0\}^{n-d}$, i.e.
we consider a ball $B(y,\nu)$ of radius $\nu$, with $y\in X_\lambda$, and we consider $X_\lambda \cap B(y,\nu)$ as the graph of a 
map $f:=f_{B(y,\nu)}$ defined on $T_yX_\lambda\cap B(y,\nu)$ as above.  We have $\mu_d(T_yX_\lambda\cap B(y,\nu)) = \nu^d$, and since $f$ is $1$-Lipschitz,
we also have $\mu_d(X_\lambda\cap B(y,\nu)) = \nu^d$.

By Lemma \ref{l.Crofton lineaire} we have the Cauchy-Crofton formula for $T_yX_\lambda\cap B(y,\nu)$ 
\begin{equation}\label{eq.6.1}
\begin{aligned}
C(n,d)\mu_d(X_\lambda\cap B(y,\nu))
&=C(n,d)\nu^d\\
&= 
C(n,d)\mu_d(T_yX_\lambda\cap B(y,\nu))\\
&=V_d(T_yX_\lambda\cap B(y,\nu)). 
\end{aligned}
\end{equation}
Using this (and $\mu_d(X_\lambda\cap B(x,\nu))\nu^{-d} = 1$), we obtain
\begin{align*}
C(n,d)\mu_d(X_\lambda)
&=
\int_{x \in X_\lambda} C(n,d) \mu_d(X_\lambda\cap B(x,\nu))\nu^{-d} 
\ \mathrm{d}x \\
&\overset{\eqref{eq.6.1}}{=} 
\nu^{-d}\int_{x \in X_\lambda} V_d(T_x X_\lambda\cap B(x,\nu)) 
\ \mathrm{d}x  \\
&=
\nu^{-d}\int_{x \in X_\lambda}
\int_{\bar{P}\in \bar{G}_n^{n-d}}
V_0(T_x X_\lambda\cap B(x,\nu)\cap \bar{P}) \ \mathrm{d}\bar{P}
\ \mathrm{d}x.  \\  
\end{align*}
In the last line, we now want to replace $T_x X_\lambda$ by $X_\lambda$ and estimate the difference this makes.
By \eqref{eq.bijection}, this changes something only for $\bar P \in \bar{\mathcal{P}}_{\nu/\lambda} = \bar{G}_n^{n-d} \setminus \bar{\mathcal{Q}}_{\nu/\lambda}$. The measure of that set is bounded by $C\nu^{d+1}/\lambda$ as observed in Step~3, and the integrands
$V_0(T_x X_\lambda\cap B(x,\nu)\cap \bar{P})$ and
$V_0(X_\lambda\cap B(x,\nu)\cap \bar{P})$ are uniformly bounded by some constant $C' \in \LL^\ZZ$ (depending only on $X$), by Remark~\ref{r.V0.bound}. Thus the error term
\[
\rho := \nu^{-d}\int_{x \in X_\lambda}
\int_{\bar{P}\in \bar{\mathcal{P}}_{\nu/\lambda,x}}
V_0(T_x X_\lambda\cap B(x,\nu)\cap \bar{P}) \ \mathrm{d}\bar{P}
\ \mathrm{d}x
\]
is bounded by $\nu^{-d}\mu_d(X_\lambda)C\nu^{d+1}/\lambda\cdot C'$, which in turn is bounded by $C''\nu$, for some $C'' \in \mathcal{C}_+(\{pt\})$ depending only on $X$ and $\lambda$ 
(again, by $\rho$ being bounded by $C''\nu$, we mean that we consider both $\rho$ and $C''\nu$ as constructible motivic functions in $\lambda$ and $\nu$,  and that we have $C''\mu - \rho \in \mathcal{C}_+(\ZZ^2)$. Note that we also have $\rho \in \mathcal{C}_+(\ZZ^2)$). 
The corresponding error term for $X_\lambda$ (instead of
$T_xX_\lambda$) satisfies the same bound; we denote it by $\rho'$, so that we have
\begin{align*}
C(n,d)\mu_d(X_\lambda) &+ \rho
= \rho' + \\
&\nu^{-d}\int_{x \in X_\lambda} 
\int_{\bar{P}\in \bar{G}_n^{n-d}}
V_0(X_\lambda \cap B(x,\nu)\cap \bar{P})\ \mathrm{d}\bar{P}
\ \mathrm{d}x.
\end{align*}
As in the proof of Lemma~\ref{l.addit}, using that for generic $\bar P$, we have 
$$V_0(X_\lambda \cap \bar P) = \mu(X_\lambda \cap \bar P)$$
and 
$$V_0(X_\lambda \cap B(x,\nu)\cap \bar P) = \mu(X_\lambda \cap B(x,\nu)\cap \bar P),$$ 
we obtain
$$\nu^{-d}\int_{x \in X_\lambda}  V_0(X_\lambda \cap B(x,\nu)\cap \bar{P}) \ \mathrm{d}x
= V_0(X_\lambda \cap \bar{P}),$$  and hence
\begin{equation}\label{eq.6.2}
\begin{aligned}
C(n,d)\mu_d(X_\lambda) + \rho
&=  \rho'  + \int_{\bar{P}\in \bar{G}_n^{n-d}}
   V_0(X_\lambda  \cap \bar{P})\ \mathrm{d}\bar{P}\\
&=   \rho' + V_d(X_\lambda).
\end{aligned}
\end{equation}
By letting $\nu$ go to $0$ and using that $\rho$ and $\rho'$ are bounded by $C''\nu$, we obtain the desired Cauchy-Crofton formula for $X_\lambda$
 \begin{equation}\label{eq.CC with lambda}
C(n,d)\mu_d(X_\lambda)=V_d(X_\lambda).
\end{equation}
Indeed, by Lemma~\ref{l.squeeze}, $\rho$ goes to $0$ for $\nu \to 0$, since the function $C''\nu$ does, and similarly for $\rho'$.
Taking the limit on both sides of \eqref{eq.6.2} (and noting that the limit is well-defined by Remark~\ref{r.fonct}) thus yields \eqref{eq.CC with lambda}.

\medskip
{\bf Step 6.}
Equation~\eqref{eq.CC with lambda} is an equation between nonnegative constructible motivic funtions in the variable $\lambda$, i.e. between functions in $\mathcal C_+(\NN)$,
so if the limit(s) exist for $\lambda \to 0$ (in the sense of Definition~\ref{d.convergence Mot}), then they are equal.

By monotone convergence (Remark \ref{r.monoton convergence}), the limit of the left hand side for $\lambda \to 0$ exists and is equal to $C(n,d)\mu_d(X)$.

To be able to apply a similar argument to the right hand side, we use once more that for generic $\bar P \in  \affgr{n-d}$, $V_0(X \cap \bar P)$ is finite and hence equal to $\mu_0(X \cap \bar P)$. We can therefore consider $V_d(X)$ as the measure of a set, and similarly for $X_\lambda$ instead of $X$. Now monotone convergence can be applied to deduce that the limit of the right hand sice of \eqref{eq.CC with lambda} is equal to $V_d(X)$, so that we finally obtain
$$C(n,d)\mu_d(X)=V_d(X),$$
as desired.
\end{proof}

\section{Motivic Vitushkin formula for metric entropy}\label{s.Motivic Vitushkin formula}

We introduce now a local version of $V_0$ in a ball $B$, which is the nonarchimedean counterpart of the local version of $V_0$ in tame real geometry as presented in introduction (see \cite[Definition 3.2]{YoCo}).
For this, we fix a $\emptyset$-definable ball $B \subseteq K^n$ and a $\emptyset$-definable set $X \subseteq K^n$. Let $B_1, \dots, B_\ell$ be as given by Theorem~\ref{t.S0} and let $S_0\subset k^m$ be the set parametrizing those balls, as in Definition~\ref{d.V0}.


\begin{notn}\label{not.ZB}
We denote by $S_0^B \subseteq k^m$ the $\emptyset$-definable subset of $S_0$ corresponding to all those $B_i$ which are contained in $B$.
\end{notn}

Note that those balls correspond to a version of Theorem~\ref{t.S0} restricted to $B$, since for any ball $B'' \subseteq B$, $X$ is not riso-trivial on $B''$ if and only if $B_i \subseteq B''$ for some $B_i$ corresponding to an element of $S_0^B$. This motivates the following relative version of $V_0$.

\begin{defn}\label{d.V_0 relative}
We set $V_0(X, B) := \mu_0(S_0^B)$, the \emph{$0$-dimensional relative motivic Vitushkin variation}.
\end{defn}

By this definition, $V_0(X, B)$ is $0$ if $X$ is $1$-riso-trivial on $B$. In the reals, this corresponds to the case where all connected components of $X \cap B$ touch the boundary of $B$.

\begin{lem}\label{l.comparison V0 and V0B}
If $X$ is not riso-trivial on $B$ or $X \cap B = \emptyset$, we have
$V_0(X, B) = V_0(X \cap B)$. Otherwise, we have $V_0(X, B) = 0$ and $V_0(X \cap B) = 1$. In particular, we always have
\[
V_0(X, B) \le V_0(X \cap B),
\]
and this inequality also holds uniformly in families, i.e. if $(X_z)_{z \in Z}$ is a $\emptyset$-definable family of sets $X_z \subseteq K^n$ which are uniformly bounded, $(B_z)_{z \in Z}$ is a $\emptyset$-definable family of balls $B_z \subseteq K^n$, then we have
\[
(V_0(X_z,B_z))_z
\le
(V_0(X_z \cap B_z))_z
\]
in $\mathcal{C}_+(Z)$.
\end{lem}

(Recall that the preorder $\le$ on $\mathcal{C}_+(Z)$ is the one from Notation~\ref{n.le}.)

\begin{proof}[Proof of Lemma~\ref{l.comparison V0 and V0B}]
Let $B_1, \dots, B_\ell$ be the balls and singletons obtained from Theorem~\ref{t.S0} for $X$ and
let $B'_1, \dots, B'_{\ell'}$ be the balls and singletons
obtained from Theorem~\ref{t.S0} for $X \cap B$.
By definition of the $B_i$ and $B'_j$, any ball $B'' \subseteq B$ appears among the $B_i$ if and any only if it appears among the $B'_j$, so the difference between $V_0(X \cap B)$ and $V_0(X, B)$ only comes from balls $B'_j$ that are not contained in $B$.
Such a $B'_j$ cannot be disjoint from $B$ (since otherwise, $X\cap B$ would be $1$-riso-trivial on that ball), so we are left to verify whether a ball $B''$ strictly containing $B$ can be among the $B'_j$. If $X \cap B = \emptyset$, then $X \cap B$ is $1$-riso-trivial on $B''$, and if $X \cap B$ is not riso-trivial on $B$, then $B''$ is not a minimal non-riso-trivial ball, so in both cases, we obtain the equality $V_0(X \cap B)= V_0(X, B)$. If we are in neither of those two cases, the smallest ball strictly containing $B$ is a minimal non-riso-trivial ball, and we obtain $ V_0(X \cap B)=1$ and $V_0(X, B) = 0$.

The family version of the inequality follows from the first part of the lemma, since the difference of the two constructible motivic functions on $Z$ is given by the indicator function of the set of those $z \in Z$ for which $X_z \cap B$ is non-empty but $X_z$ is $1$-riso-trivial on $B_z$.
\end{proof}

For $0 < i \le n$, we then define $V_i(X, B)$, the \emph{$i$th  relative motivic Vitushkin variation}, analogously to Definition~\ref{d.Vi}.

\begin{defn}\label{d.V_i relative}
For $X,B$ as above and $1 \le i \le n$, set
$$
V_i(X,B) :=  \int_{\bar P \in \affgr{n-i},\bar P \cap B \ne \emptyset} V_0(X \cap \bar P,B \cap \bar P) \ \mathrm{d} \bar P.
$$
\end{defn}

The following Proposition 
is an analogue of \cite[Prop. 3.3]{YoCo}, which asserts that for any 
disjoint sets $B_1,\ldots, B_\ell\subseteq \RR^n$, for any set $X\subseteq \RR^N$ and for any $i\in \{0,\ldots, n\}$,  one has
$\sum_{j=1}^\ell V_i(X,B_j) \le V_i(X)   $.

\begin{prop}
\label{p.bound of integral of ViB}
For any $i \le n$, any bounded $\emptyset$-definable set $X\subseteq K^n$ and any  $\lambda \in  \LL^\ZZ$, 
\begin{equation}\label{eq.int ViB}
\lambda^{-n} \int_{x\in K^n} V_i(X,  B(x,\lambda)) \ \mathrm{d}x \le V_i(X).
\end{equation}
The same inequality also holds in families, that is, given a $\emptyset$-definable set $Z$, a $\emptyset$-definable function $Z  \to \LL^\ZZ, z \mapsto \lambda_z$ and a $\emptyset$-definable family
$(X_z)_{z \in Z}$ of sets $X_z\subseteq K^n$,  with $X_z$ uniformly bounded, we have 
\begin{equation}\label{eq.int ViB relatif}
\left( \lambda_z^{-n} \int_{x\in K^n} V_i(X_z,  B(x,{\lambda_z})) \ \mathrm{d} x\right)_z
 \le (V_i(X_z))_z
\end{equation}
in $\mathcal{C}_+(Z)$.
\end{prop}

\begin{proof}
For simplicity, we formulate the entire proof in the non-family case; the family-version works in exactly the same way.

First consider the case $i = 0$. The case $n = 0$ is trivial, so suppose that $n \ge 1$.

Let $B_1, \dots, B_\ell$ be the balls obtained from Theorem~\ref{t.S0} for $X$, let $S_0$ be the set parametrizing those balls as in Definition~\ref{d.V0} (so that $V_0(X) = \mu_0(S_0)$), and denote by
$S_0^{B(x,\lambda)}$ the subset of $S_0$ corresponding to those balls $B_i$ contributing to $V_0(X, B(x,\lambda))$, i.e. satisfying $B_i \subseteq B(x,\lambda)$. Those sets $S_0^{B(x,\lambda)}$ are either disjoint or equal, and their union $T_0 := \bigcup_x S_0^{B(x,\lambda)}$ corresponds exactly to those balls $B_i$ which have radius at most $\lambda$.

The left hand side of \eqref{eq.int ViB} can be rewritten as
\begin{align*}
\lambda^{-n} \int_{x\in K^n} V_0(X,  B(x,\lambda)) \ \mathrm{d}x
&=
\lambda^{-n} \int_{x\in K^n} \mu_0(S_{0}^{B(x,\lambda)}) \ \mathrm{d} x
\\
&=
\lambda^{-n}\sum_{s \in T_{0}} \mu(\{x ; s \in S_{0}^{B(x,\lambda)}\})\\
&=\mu_0(T_0),
\end{align*}
where the sum sign stands for the 0-dimensional motivic integral, and where in the last equality, we use that for each $s \in T_0$, the set in the sum is a ball of radius $\lambda$ and hence has measure $\lambda^n$ (which cancels with the $\lambda^{-n}$).
This proves the proposition in the case $i = 0$, 
applying Remark \ref{r.subset and inequality} to the inclusion 
$S_0^{B(x,\lambda)}\subset S_0$.

For the case $i \ge 1$, we apply the case $i=0$ to the $\emptyset$-definable families 
$(X_z \cap \bar P)_{z\in Z, \bar{P}\in \affgr{n-i}}$ and $(\lambda_z)_{z\in Z}$. We obtain 
$$
\int_{x\in K^n} V_0(X_z \cap \bar{P},  B(x,\lambda_z)) \ \mathrm{d} x \le 
\lambda_z^n   V_0(X_z\cap \bar{P})
.$$ 
By preservation of the partial order $\le$ under integration, we obtain
$$
\int_{\bar{P}\in \affgr{n-i}} \int_{x\in K^n} V_0(X_z\cap \bar{P},  B(x,\lambda_z)) \ \mathrm{d} x \ \mathrm{d}\bar{P} 
\le
\int_{\bar{P}\in \affgr{n-i}}
 \lambda_z^n  V_0(X_z\cap \bar{P})
 \ \mathrm{d}\bar{P} 
.$$ 
But by Fubini theorem the left hand side of the above inequality is
$$ 
\int_{x\in K^n}
\int_{\bar{P}\in \affgr{n-i}}  
V_0(X_z\cap \bar{P},  B(x,\lambda_z)) \ \mathrm{d}\bar{P} \ \mathrm{d} x 
=  \int_{x\in K^n}  
V_i(X_z,  B(x,\lambda_z)) \ \mathrm{d} \bar{P} \ \mathrm{d} x,
$$
and the right hand side is by definition $\lambda_z^n V_i(X_z)$. 
\end{proof}

The following is an analogue of inequality \eqref{eq.EntRelPol} in the introduction (see \cite[Thm. 3.4]{YoCo}). This is the first place where we cannot work anymore with the partial preorder $\le$ on constructible motivic functions, but need to work with the preorder $\preccurlyeq$ instead. The reason is the following: so far we could use inclusions and Remark \ref{r.subset and inequality} to produce bounds with respect to the partial preorder $\leq$. But in the proof of Theorem \ref{t.bound on sum of variations}, we have to compare the measure $\mu_0(S_0^B)$, as an element 
of $\mathcal{C}_+(\{pt\})$, to the constant function $1$ on the point, and in this situation, where inclusions provide no help, one has to replace them by surjections (from $S_0^B$ to $\{pt\}$).

\begin{thm}\label{t.bound on sum of variations}
Let $n \in \NN$ be given.
There exists an invertible element $C(n) \in C_+(\{pt\})$ (depending only on $n$) such that
for any $\emptyset$-definable set $X\subseteq K^n$ and any  $\lambda \in  \LL^\ZZ$,
assuming that $B$ is a $\emptyset$-definable ball of radius $\lambda$ such that $X \cap B \ne \emptyset$, we have 
$$
\sum_{i = 0}^n  \lambda^{-i} V_i(X, B) \succcurlyeq C(n).
$$
The same inequality also holds in families, that is, given a $\emptyset$-definable measurable set $Z\subseteq K^m$,
 a $\emptyset$-definable function $Z  \to \LL^\ZZ, z \mapsto \lambda_z$,
two $\emptyset$-definable families
$(X_z)_{z \in Z}$ of sets $X_z\subseteq K^n$,  with $X_z$ uniformly bounded, and 
$B_z$ a $\emptyset$-definable ball of radius $\lambda_z$ such that $X \cap B \ne \emptyset$,
we have 
$$\left( \sum_{i = 0}^n  \lambda_z^{-i} V_i(X_z, B_z) \right)_z \succcurlyeq C(n)$$
in $\mathrm{I}\mathcal{C}(Z)$.
\end{thm}

Before starting with the main proof, we compute some measure which will make up the constants $C(n)$. Note that the precise value depends on the normalization of the measure on $ \vecgr{n-d}$ we fixed in Notation~\ref{n.normalisation measure}.

\begin{lem}\label{l.G}
Let $d \le n$ be given, set $V := \{0\}^{n-d}\times k^d \subset k^n$, and let $G \subset  \vecgr{n-d}$ be the set of those vector spaces $U \subseteq K^n$ satisfying $\res(U) \cap V = \{0\}$. Then the measure of $G$ is equal to
$$\mu(G) = \prod_{i=1}^d (1-\LL^{-i}) \cdot \prod_{i=1}^{n-d} (1-\LL^{-i}).$$
\end{lem}

\begin{proof}
By our definition of the measure on $ \vecgr{n-d}$ (see Notation 
\ref{n.mesure grassmannienne}), we have $\mu(G)=\mu_{\GL_n}(H)$, where $H = \{g\in \GL_n(\mathcal{O}_K); g\cdot P_0\in G\}$,
and where we use $P_0=K^d\times 0$ as fixed in Remark~\ref{r.identifying Ki and Kn/P}.
Then $\res(g\cdot P_0)$ is transverse to $V$ if and only if the top left $d\times d$ minor of $g$ is invertible. 
This implies
$$H=\begin{pmatrix}
1&0\\Mat_{d\times(n-d)}(\OK)&1
\end{pmatrix}\cdot
\begin{pmatrix}
\GL_d(\OK)&Mat_{(n-d)\times d}(\OK)\\0&\GL_{n-d}(\OK)
\end{pmatrix}.
$$
One verifies that the map
\begin{align*}
&Mat_{d\times(n-d)}(\OK) \times \GL_d(\OK) \times Mat_{(n-d)\times d}(\OK)\times \GL_{n-d}(\OK) \to  H \hfill
\\
&(C,A,B,D) \mapsto \begin{pmatrix}
1&0\\C&1
\end{pmatrix}\cdot
\begin{pmatrix}
A&B\\0&D
\end{pmatrix}
\end{align*}
is measure-preserving, so we have
\begin{align*}
    \mu(H) =& \mu(Mat_{d\times(n-d)}(\OK)) \cdot \mu(\GL_d(\OK))\\
  &  \cdot\mu(Mat_{(n-d)\times d}(\OK))\cdot \mu( \GL_{n-d}(\OK))\\
=& 1\cdot [\GL_d(k)]\cdot\LL^{-d^2}\cdot 1
\cdot  [\GL_{n-d}(k)]\cdot\LL^{-(n-d)^2}.
\end{align*}
The lemma now follows using
\[[\GL_d(k)] = \prod_{i=0}^{d-1} (\LL^d-\LL^i).
\qedhere
\]
\end{proof}

\private{Probably the following should be true: If we use the canonical measure on the grassmanians (satisfying $\mu( \vecgr{n-d}) = [ \vecgr{n-d}]\cdot\LL^{-dim}$), then the map $\mathrm{Mat}_{d\times(n-d)}(\OK) \to  \vecgr{n-d}$ sending a matrix to its graph should be measure-preserving, and its image is the set $G$ from the proof. This shows that $G$ has measure $1$, and so we obtain that the constant in \ref{t.bound on sum of variations} can be taken to be 1.}

\begin{proof}[Proof of Theorem~\ref{t.bound on sum of variations}]
We mostly formulate the proof in the absolute version; for the family version we just indicate the necessary changes.

By $i$-homogeneity of $V_i(\cdot, B)$ (see Remark \ref{r.Vi}\eqref{item.homogeneity}), we may assume $\lambda = 1$.
Let $V:=\rtsp_B(X) \subseteq k^n$ be the riso-triviality space of $X$ on $B$, as in Definition~\ref{d.risotriv} and set $d:=\dim V$. To prove the theorem, we set $C(n) := (1-\LL^{-1})^n$, and we will prove that one single summand is bigger or equal to $C(n)$, namely $V_d(X, B) \succcurlyeq C(n)$.

In the case $d = 0$, $X$ is not riso-trivial on $B$,
so the set $S_{0}^B$ appearing in Definition \ref{d.V_0 relative} of $V_0(X, B)$ is non-empty. In particular there exists a surjection 
$S_{0}^B \longtwoheadrightarrow \{0\}$, so we have
$$V_0(X, B) = \mu_0(S_{0}^B) \succcurlyeq 1 \succcurlyeq C(n).$$
In the family version, let us write $S_{0,z}^B$ for the set corresponding to $X_z$. We then have a family of surjections $S_{0,z}^B \longtwoheadrightarrow \{0\}$;
to deduce $(\mu_0(S_{0,z}^B))_z \succcurlyeq 1$, we use Lemma \ref{l.<surj}.

In the case $d \ge 1$, let $G \subseteq \gr{n-d}$ consist of those vector spaces $U \subseteq K^n$ satisfying
$\res(U) \cap V = \{0\}$.
Recall that $V$ is $\emptyset$-definable by Theorem~\ref{t.defble}, so we may as well assume that $V = \{0\}^{n-d} \times K^d$, so that $G$ becomes the same $G$ as in Lemma~\ref{l.G}, which implies
$\mu(G) \ge C(n)$.

Let $\bar G \subseteq \affgr{n - d}$ be the definable set of affine spaces of $K^n$ with direction in $G$, and $\bar{G}(B)$ the set of affine spaces $\bar{P}\in \bar{G}$ such that $\bar{P}\cap B\not=\emptyset$. The measure of that set is equal to the measure of $G$, so we also have $\mu(\bar G(B)) \ge C(n)$.

Now consider $\bar P \in \bar G(B)$. By Proposition~\ref{p.fib},
$X \cap B \cap \bar P$ is not riso-trivial on $B \cap \bar P$, so we have that $V_0(X \cap \bar P, B \cap \bar P) \succcurlyeq 1$, as in the case $d = 0$, uniformly in $\bar{P}$ by Lemma \ref{l.<surj} (and also uniformly in $z$, in the family version).
Therefore, integrating over $\bar G$ (and using Lemma~\ref{l.<int}), we obtain
\begin{equation}\label{eq.ge 1}
V_d(X, B) \succcurlyeq \mu(\bar{G}(B))\succcurlyeq C(n),
\end{equation}
as desired.
\end{proof}

We now define, for $\lambda=\LL^r \in \LL^\ZZ$, the motivic
$\lambda$-entropy $M(X,\lambda)$ of a $\emptyset$-definable subset $X$ of $K^n$. This
is the nonarchimedean counterpart of the real metric entropy of the introduction, since this is the motivic measure of balls of radius $\lambda$ one needs to cover $X$, but this is also the measure of the tubular neighbourhood $T_r(X)=\{x\in K^n; \exists y\in X, \val(y-x)\ge r\}$ of $X$.

\begin{defn}[Motivic entropy]\label{d.Motivic entropy}
Let $\lambda \in \LL^\ZZ$ and $X$ a $\emptyset$-definable subset of $K^n$.
The \emph{ $\lambda$-entropy $M(X, \lambda)$  of $X$} is the measure 
$$
M(X, \lambda) :=\lambda^{-n} \mu(\{x ; B(x,\lambda) \cap X \ne \emptyset\})
$$
\end{defn}

The following Theorem is the nonarchimedean analogue of inequality 
\eqref{eq.EntPol} given in the introduction
(see \cite[Thm 3.5]{YoCo}) and the main result of this paper.
\begin{thm}[Ivanov's bound for motivic entropy]
\label{t.ent}
Let $X$ be a boun\-ded $\emptyset$-definable subset of $K^n$.  
For any $\lambda \in \LL^\ZZ$ one has
$$
M(X, \lambda)C(n) \preccurlyeq \sum_{i=0}^n \lambda^{-i} V_i(X),
$$
where $C(n)\in \mathcal C_+(\{pt\})$ is the same invertible constant as in Theorem~\ref{t.bound on sum of variations}.

The same inequality also holds in families, that is, given a $\emptyset$-definable measurable set $Z\subseteq K^m$, 
 a $\emptyset$-definable function $Z  \to \LL^\ZZ, z \mapsto \lambda_z$,
a $\emptyset$-definable family
$(X_z)_{z \in Z}$ of sets $X_z\subseteq K^n$,  with $X_z$ uniformly bounded, 
we have 
$$
\left( M(X_z, \lambda_z) \right)_zC(n) \preccurlyeq  \left( \sum_{i=0}^n \lambda^{-i} V_i(X)\right)_z.
$$
in $\mathrm{I}\mathcal{C}(Z)$.
\end{thm}

\begin{proof}
We  have by Proposition \ref{p.bound of integral of ViB}
\begin{align*}
\sum_{i=0}^n \lambda^{-i} V_i(X) & \ge \sum_{i=0}^n  
\lambda^{-i-n} \int_{x\in K^n} V_i(X,  B(x,\lambda)) \ \mathrm{d}x \\
& =   \lambda^{-n} \int_{x\in K^n}\sum_{i=0}^n \lambda^{-i}  V_i(X,  B(x,\lambda)) \ \mathrm{d}x.
\end{align*}

Denoting  $Y := \{x ; B(x,\lambda) \cap X \ne \emptyset\}$, 
for $x \in Y$, Theorem~\ref{t.bound on sum of variations} yields that the integrand of the last integral above  is bounded from below by $C(n)$, so we obtain

\[\sum_{i=0}^n \lambda^{-i} V_i(X)
  \succcurlyeq  \lambda^{-n} \mu(Y)C(n)  = M(X,\lambda)C(n).
  \qedhere
\]
\end{proof}

\section{Appendix: Tensor products of semimodules}
\label{s.tensor}

 In the definition of motivic integration in \cite{CL.mot}, a notion of tensor product of two semirings over a common subsemiring appears. From a category theoretical point of view, what is meant is simply the pushout in the category of commutative semirings. In this appendix, we describe this pushout more explicitly, as needed in the proof of Lemma~\ref{l.trivial-writing} below.
Note that   a 
definition of tensor product over semirings is given in \cite{Gol.semiRing} (which sends back to \cite{Taka.semiTensor}), but that definition has a slightly different universal property than the one we want; in particular, it does not yield the pushout of semirings.

We define a semiring $S$ and an $S$-semimodule $M$ as in \cite{Gol.semiRing}; in particular, $S$ contains neutral elements $0$ and $1$, and we require that for every $m \in M$, we have $0 \cdot m = 0$ and $1 \cdot m = m$. All our semirings are commutative.

The free $S$-semimodule over a set $A$ is the set of maps from $A$ to $S$ with finite support. Following \cite{Gol.semiRing}, we denote it by $S^{(A)}$.

Let $M_1$ and $M_2$ be $S$-semimodules and set $L := S^{(M_1 \times M_2)}$. Given $m_1 \in M_1, m_2 \in M_2$, we write $[m_1,m_2] \in L$ for the natural image of $(m_1,m_2)$ in $L$, i.e. $[m_1,m_2]$ is the map taking the value $1$ at $(m_1,m_2)$ and $0$ everywhere else.

Let $P_0 \subset L \times L$ be the set of pairs of the form
\begin{enumerate}
    \item 
$(s[m_1, m_2], [sm_1, m_2])$
\item $(s[m_1, m_2], [m_1, sm_2])$
\item $([m_1+m_1',m_2], [m_1,m_2]+[m_1',m_2])$;
\item $([m_1,m_2+m_2'], [m_1,m_2]+[m_1,m_2'])$
\end{enumerate}
for $m_i, m'_i \in M_i$ and $s \in S$, set $P_1 := P_0 \cup \{(a,b); (b,a) \in P_0\}$, and let $P_2$ be the $S$-sub-semimodule of $L\times L$ generated by $P_1$, i.e. consisting of sums of the form
\begin{equation}\label{eq.sbb}
\sum_{j=1}^m s_j(b_j, b'_j)
\end{equation}
for $(b_j,b'_j) \in P_1$ and $s_j\in S$.
Note that $P_2$ can be considered as a symmetric and reflexive relation on $L$. We define $\simeq$ to be the transitive closure of that relation, that is to say, for $a, a' \in L$,
we have $a \simeq a'$ if and only if we have
\begin{equation}\label{eq.aaa}
a = a_0, \dots, a_n = a',
\end{equation} such that for each $i$, we can write
$a_i = \sum_j s_jb_j$, $a_{i+1} = \sum_j s_jb'_j$ for some $s_j \in S$ and some $(b_j, b'_j) \in P_1$.
Clearly, $\simeq$ is an equivalence relation on $L$, and $L/\mathord{\simeq}$ carries a natural $S$-semimodule structure, since $a_1 \simeq b_1$ and $a_2 \simeq b_2$ implies $a_1+a_2 \simeq b_1+b_2$ and $sa_1 \simeq sb_1$, for $a_i, b_i \in L$ and $s \in S$.

\begin{defn}\label{d.tensor}
Let $M_1$ and $M_2$ be $S$-semimodules over a commutative semiring $S$. Then, using the above notation, we define the tensor product of $M_1$ and $M_2$ to be the $S$-semimodule $$M_1 \otimes_{S} M_2 := L/\mathord{\simeq}.$$
The image of $[m_1,m_2] \in L$ in $M_1 \otimes_{S} M_2$ is denoted by $m_1\otimes m_2$.
\end{defn}

\begin{rem}
In contrast to \cite[Proposition~(16.12)]{Gol.semiRing}, with our definition, the tensor product $M_1 \otimes_{S} M_2$ is not necessarily cancellative, and accordingly, it satisfies a slightly different universal property than \cite[Proposition~(16.14)]{Gol.semiRing}. Also, our definition is slightly different because we assume $S$ to be commutative.
\end{rem}

\begin{rem}
If $M_1, M_2$ and $S$ are commutative semirings coming with (semiring) morphisms $S \to M_i$ (which turn $M_i$ into $S$-semimodules), then $M_1 \otimes_S M_2$ is the pushout of $M_1$ and $M_2$ over $S$ in the category of commutative semirings. Indeed, 
$M_1 \otimes_S M_2$ carries a natural multiplication (satisfying $(m_1 \otimes m_2)\cdot (m'_1 \otimes m'_2) = m_1m'_1 \otimes m_2m'_2$, it comes with morphisms $M_i \to M_1 \otimes_S M_2$ which agree on $S$ (namely, $m_1 \mapsto m_1 \otimes 1$ and $m_2 \mapsto 1 \otimes m_2$ for $m_i \in M_i$), and it satisfies the universal property of the pushout: for any other commutative semiring $T$ coming with morphisms $M_1 \to T, M_2 \to T$ which agree on $S$, there exists a unique morphism $M_1 \otimes_S M_2 \to T$ making everything commute. We leave it to the reader to verify all those claims.
\end{rem}

The following lemma is the reason why we go through all this hassle. We are interested in the following property ($*$) of a pair of $S$-semimodules $U \subset M$ (where $U$ is an $S$-sub-semimodule of $M$):
\begin{itemize}
    \item[($*$)] if $m,m'\in M$ and $s \in S \setminus \{0\}$ satisfy $sm+m' \in U$, then $m,m' \in U$.
\end{itemize}

\begin{lem}\label{l.trivial-writing}
Let $S$ be a (commutative) semiring, let $M_1$ and $M_2$ be $S$-semimodules,
let $U_i \subset M_i$ be $S$-sub-semimodules, for $i=1,2$,
and let $U \subset M_1 \otimes_{S} M_2$ be the $S$-sub-semimodule generated by elements of the form 
$u_1 \otimes m_2$ and $m_1 \otimes u_2$, for $u_i \in U_i$ and $m_i \in M_i$.
If $\{0\} \subset S$ satisfies ($*$) and both $U_1 \subset M_1$ and $U_2 \subset M_2$ satisfy ($*$), then $U \subset M_1 \otimes_{S} M_2$ also satisfies ($*$).
\end{lem}

\begin{proof}
Let $L$, $P_i$, $\simeq$ be defined as above Definition~\ref{d.tensor}, and let $\tilde U$ be the $S$-sub-semimodule of $L$ generated by elements of the form $[u_1, m_2]$ and $[m_1, u_2]$, for $u_i \in U_i$ and $m_i \in M_i$.
Using that $\{0\} \subset S$ satisfies ($*$), one deduces that also $\tilde U \subset L$ satisfies ($*$).

To prove the lemma, it suffices to show that $B'$ is closed under $\simeq$, i.e. $a \in \tilde U$ and $a \simeq a'$ implies $a' \in \tilde U$. Indeed, given $m,m',s$ as in ($*$) for $U\subset M_1 \otimes_{S} M_2$, choose representatives $\tilde m, \tilde m' \in L$. The assumption $sm + m' \in U$ implies that $sm + m'$ has a representative $a' \in \tilde U$.
We have $s\tilde m + \tilde m' \simeq a'$, so using that $\tilde U$ is closed under $\simeq$, we deduce $s\tilde m + \tilde m' \in \tilde U$. Now ($*$) for $\tilde U \subset L$ implies $\tilde m, \tilde m' \in \tilde U$ and hence $m, m' \in U$.

To show that $\tilde U$ is closed under $\simeq$, it suffices to verify this for individual steps in the sequence \eqref{eq.aaa}, i.e. we need to show that for $s_j, b_j, b'_j$ as in \eqref{eq.sbb}, $\sum_j s_jb_j \in \tilde U$ implies 
$\sum_j s_jb'_j \in \tilde U$, so suppose that $\sum_j s_jb_j \in \tilde U$.
Since $\tilde U \subset L$ satisfies ($*$), this implies that each $b_j$ lies in $\tilde U$.
It therefore suffices to show that, for each pair $(b,b') \in P_0$, we have $b \in \tilde U$ if and only if $b' \in \tilde U$.
This is indeed the case in all cases of the definition of $P_0$:

Case (1): If $s = 0$, then both sides lie in $\tilde U$, so suppose now that $s \ne 0$. Then we have $s[m_1, m_2] \in \tilde U \iff [m_1,m_2] \in \tilde U \iff [sm_1,m_2] \in \tilde U$, where the first equivalence follows from ($*$) for $U \subset L$, and the second one is trivial if $m_2 \in U_2$ and follows from ($*$) for $U_1 \subset M_1$ if $m_2 \notin U_2$.

Case (2): Analogous.

Case (3): This is similar:
By ($*$) for $U \subset L$, we have $[m_1,m_2]+[m_1',m_2] \in \tilde U$ if and only if both, $[m_1,m_2]$ and $[m_1',m_2]$ lie in $\tilde U$,
and this is equivalent to
$[m_1+m_1',m_2] \in \tilde U$, either trivially (if $m_2 \in U_2$) or using ($*$) for $U_1 \subset M_1$ (if $m_2 \notin U_2$).

Case (4): Analogous.
\end{proof}

\bibliographystyle{siam}
\bibliography{references}

\vfill
\end{document}